\documentclass[a4paper,11pt]{article}     % `article' for A4 paper
\usepackage{array}                        % array extensions

\usepackage{amsmath,amssymb}        % AMS-LaTeX with CD and symbols
\usepackage{amsthm}
\usepackage{float}  
\input xy
\xyoption{all}
\usepackage[dvips]{color}  

\usepackage[dvips]{graphicx}

\usepackage{rotating}

\DeclareMathOperator{\codim}{codim}
\DeclareMathOperator{\Frob}{Frob}
\DeclareMathOperator{\trace}{trace}

\newtheorem{theorem}{Theorem}[section]
\newtheorem{lemma}[theorem]{Lemma}
\newtheorem{corollary}[theorem]{Corollary}
\newtheorem{proposition}[theorem]{Proposition}

\theoremstyle{remark}
\newtheorem{remark}[theorem]{\bf Remark}

\newtheorem{definition}[theorem]{\bf Definition}

\newcommand{\AAA}{{\mathbb{A}}}
\newcommand{\CC}{{\mathbb{C}}}

\newcommand{\FF}{{\mathbb{F}}}
\newcommand{\PP}{{\mathbb{P}}}
\newcommand{\QQ}{{\mathbb{Q}}}
\newcommand{\RR}{{\mathbb{R}}}
\newcommand{\ZZ}{{\mathbb{Z}}}

\newcommand{\ba}{{\boldsymbol{a}}}
\newcommand{\bb}{{\boldsymbol{b}}}
\newcommand{\bc}{{\boldsymbol{c}}}

\newcommand{\eps}{{\varepsilon}}
\newcommand{\ra}{{\rightarrow}}

\newcommand{\comment}[1]{}

\newcommand{\rl}[1]{\rule{0ex}{#1ex}}

\newcommand{\noteone}{
\left\{\!\!
\begin{array}{l}
{\mathcal F}_1:
\text{If }\sum\sqrt{a_i} = 0,\\
X_\ba \text{ has }\ge 31\\
\text{singularities,}\\
\text{and }h^{12}(\widetilde X_\ba)\le 4
\end{array}
\right.
}

\newcommand{\notetwo}{
\left\{\!\!
\begin{array}{l}
\text{If }\ba=(1\!:\!1\!:\!1\!:\!1\!:\!1\!:\!t),\\
t\not=0,1,9,\text{ then }\\
\dim W_\ba=8, \text{if also }\\
t\not=25, \text{ then }h^{12}=5
\end{array}
\right.
}

\newcommand{\notefive}{
\begin{array}{l}
\text{If }\ba=(1\!:\!1\!:\!1\!:\!r\!:\!r\!:\!r),\\
r\not=0,1,9,\text{ then }\\
\dim W_\ba=8, h^{12}=5.\\
\text{if  }r=9, \text{ then }\\
\dim W_\ba=4, h^{12}=2
\end{array}
}

\newcommand{\notethree}{
\left.
\begin{array}{r}
{\mathcal F}_{10}:
\text{If }\ba=(1\!:\!1\!:\!b^2\!:\!b^2\!:\!(b+1)^2\!:\!(b+1)^2),\\
X_\ba\text{ has }\ge 35\text{ singularities }\\
\text{and }h^{12}(\widetilde X_\ba)\le 1
\end{array}
\!\!
\right\}
}

\newcommand{\notefour}{
\left\{
\begin{array}{l}
{\mathcal F}_8:\text{If }\ba=(1\!:\!1\!:\!1\!:\!b^2\!:\!4\!:\!(b+1)^2),\\
b\not=0,1,\;X_\ba\text{ has }\ge 34\\
\text{singularities, and }h^{12}(\widetilde X_\ba)\le 1
\end{array}
\!\!
\right.
}

\begin{document}
\title{On modularity of rigid and nonrigid
Calabi-Yau varieties associated to the root lattice $A_4$}
\author{Klaus Hulek and Helena Verrill}
\date{}

\maketitle

\begin{abstract}
We prove the modularity of
four rigid and three nonrigid Calabi-Yau threefolds
assciated with the $A_4$ root lattice.
\end{abstract}

\section{Introduction}

In this paper we investigate the geometry and arithmetic of a family
of Calabi-Yau threefolds $X_\ba$, $\ba=(a_1:\cdots:a_6)\in\PP^5$,
birational to the projective hypersurface in $T:=\PP^4
\setminus\{X_1\dots X_5=0\}$ given by 
$$
X_\ba\cap T:\;\;
(X_1+\ldots +X_5)\left(\frac {a_1}{X_1}+\ldots +\frac {a_5}{X_5}\right)=a_6.
$$
Our motivation
is to find further examples of modular Calabi-Yau 
varieties, i.e., of Calabi-Yau varieties which are defined over the 
rationals and whose $L$-series can be described in terms of modular
forms.
This 
is motivated by the
Fontaine-Mazur conjecture on the modularity of two dimensional
$\ell$-adic Galois representations coming from geometry 
\cite{FM},
which is a generalization of the 
Taniyama-Shimura-Weil conjecture on the modularity of elliptic curves,
proved by 
%Breuil, Conrad, Diamond, Taylor, and Wiles,
Wiles et al, 
\cite{Wi}, \cite{BCDT}.
More precisely, Fontaine and Mazur \cite{FM}
define the notions of a \lq\lq geometric Galois representation\rq\rq,
and a \lq\lq Galois representation coming from geometry\rq\rq.
They conjecture that geometric Galois representations 
are precisely the Galois representations coming from geometry
(such as the ones we consider) (\cite[conjecture 1]{FM}),
and combining this with classical conjectures (see e.g. \cite{serre2}) leads them to the 
conjecture that two-dimensional irreducible geometric Galois representations
are modular up to a Tate twist (\cite[conjecture 3c]{FM}).
From these two conjectures, one obtains the conjecture that
two dimensional irreducible Galois representations coming from geometry 
are modular up to a Tate twist.  We will use the term \lq\lq modular \rq\rq to mean 
\lq\lq modular up to a Tate twist \rq\rq, and even to denote direct sums of modular 
Galois representations.

Rigid Calabi-Yau threefolds (defined over $\QQ$) are expected to be modular,
since they have $2$-dimensional middle cohomology.
One expects the $L$-series
of the Galois action on the middle $\ell$-adic cohomology
to be the Mellin transform of a weight $4$ elliptic modular form.
Although recently Dieulefait and Manoharmayum \cite{DM} proved
that a rigid Calabi-Yau threefold is modular provided
it has good reduction at $3$ and $7$, 
or at $5$ and another suitable prime (in fact
most of our examples have bad reduction
at $3$ and $5$, so this result does not apply, and in general does not
determine the exact modular form),
it is still the case that relatively few examples of modular Calabi-Yau
threefolds are explicitly known.
Most currently known examples are given in
Yui's survey articles \cite{yui1}, \cite{yui2}.
Other recent examples are given by \cite{CM}.

Modularity has been also conjectured
for certain nonrigid examples, e.g., \cite{CM}.
What is new in this paper is the proof of modularity for
several nonrigid examples.
Note that we mean modularity in the sense of 
Fontaine-Mazur, i.e.,
the semisimplification of the 
Galois representation is a sum of $2$ dimensional pieces.
There are  few examples of other kinds of modularity 
of nonrigid Calabi-Yau threefolds 
known. 
Consani and Scholten \cite{consani-scholten}
consider an example corresponding to a Hilbert modular form
for which they provide evidence for the modularity and
Livn\'e and Yui \cite{livne-yui} very recently
gave some cases
involving weight $2$ and $3$ forms. Their
examples and techniques are quite different from ours.

We shall study a certain
$5$-dimensional family $X_\ba$, $\ba\in\PP^5$,
of (singular) Calabi-Yau threefolds,
associated to the root lattice $A_4$,
by means of
Batyrev's construction  \cite{Ba} of Calabi-Yau varieties as toric
hypersurfaces. 

The rigid cases are $X_1$, $X_9$,
$X_{(1:1:1:1:4:4)}$, and $X_{(1:1:1:4:4:9)}$,
where $X_t:=X_{(1:1:1:1:1:t)}$.  These have
$40$, $35$, $37$ and $35$ nodes respectively,
and their (big) resolutions have $2$-dimensional middle cohomolgy.
We will show that the $L$-series is 
%(up to the factors associated to bad primes)
the Mellin transform of a modular form of weight $4$, and level 
$6$, $6$, $12$ and $60$ respectively.  
The first few terms of their $q$-expansions are 
\begin{eqnarray}
f_6 &=&   q - 2q^2 - 3q^3 + 4q^4 + 6q^5 + 6q^6 - 16q^7 - 8q^8 +
\dots,
\label{eqn:f6}
\\
f_{12} &=&q + 3q^3 - 18q^5 + 8q^7 + 9q^9 + 36q^{11} - 10q^{13} +
\dots,
\label{eqn:f12}
\\
f_{60} &=& q - 3q^3 - 5q^5 - 28q^7 + 9q^9 - 24q^{11} - 70q^{13} + 
15q^{15} +\dots,
\label{eqn:f60}
\end{eqnarray}
where $f_N$ has level $N$.
Although the middle cohomology of $X_1$ and $X_9$ have the same L-series,
we will see that they are not birational to each other, though
by a conjecture of Tate, one expects a correspondence
between them.
The nonrigid examples we consider are
  $X_{25}$, $X_{(1:1:1:9:9:9)}$ and $X_{(1:1:4:4:4:16)}$.
  In these cases we show that $L$-series
  of the middle cohomology of the big resolutions are
  \begin{equation}
  \label{eq:L-series_sum}
  L(f_{30},s){L(g_{30},s-1)}^4,\;\;
  L(f_{30}',s){L(g_{30},s-1)}^2,\;\;
  L(f_{90},s)L{(g_{30},s-1)},
  \end{equation}
  respectively,
  where $L(h)$ denotes the Mellin transform of the function $h$, and
  the functions $g_{30}, f_{30}, f_{30}'$ and $f_{90}$ are
  cuspidal Hecke eigen newforms, $g_{30}$ having weight $2$, the others
  weight $4$, and $f_{90}$ having level $90$, the others level $30$.
(The  level $30$ has, to our knowledge, not previously
appeared in examples of this kind.)
The first few terms of the $q$-expansions are
\begin{eqnarray}
g_{30}(q) &=& q - q^2 + q^3 + q^4 - q^5 - q^6 - 4q^7 - q^8 + q^9 + q^{10} 
+ q^{12} +
\dots
\label{eqn:g30}
\\
f_{30}(q) &=& q - 2q^2 + 3q^3 + 4q^4 + 5q^5 - 6q^6 + 32q^7 - 8q^8 + 9q^9 +\dots
\label{eqn:f30}
\\
f_{30}'(q)&=& q + 2q^2 + 3q^3 + 4q^4 - 5q^5 + 6q^6 - 4q^7 + 8q^8 + 9q^9  +\dots
\label{eqn:f30p}
\\
f_{90}(q) &=& q - 2q^2 + 4q^4 - 5q^5 - 4q^7 - 8q^8 + 10q^{10} - 12q^{11}
 +\dots
\label{eqn:f90}
\end{eqnarray}
Given expression (\ref{eq:L-series_sum}), one would expect,
by the Tate conjecture, that there is a
geometric reason for the occurrence of the weight $2$ modular form
$g_{30}$, which is the Mellin transform of the L-series of a certain
elliptic curve.  We will see that this is indeed the case.

In both rigid and nonrigid cases, we use the powerful theorem
due to Faltings, Serre and Livn\'e  \cite{Livne},
which permits one to determine 
$2$-dimensional Galois representations from a finite set of data.
In practice this means counting the number of points modulo $p$ for a
given
finite number of primes, a task which can be done easily by computer. 

In \S 2 we consider the toric geometry set up.  In \S 3 we discuss the
resolution of singularities of the singular subfamily $X_\ba$.
In \S 4 we show that $X_\ba$ is birational to a fibre product of
families of elliptic curves, which allows us to apply results of
Schoen \cite{schoen}.  In \S 5 we study a certain elliptic surfaces
contained in $X_\ba$, and in \S 6 we count points and apply Livn\'e's
method to determine the L-series of the $7$ cases of $X_\ba$ mentioned above.

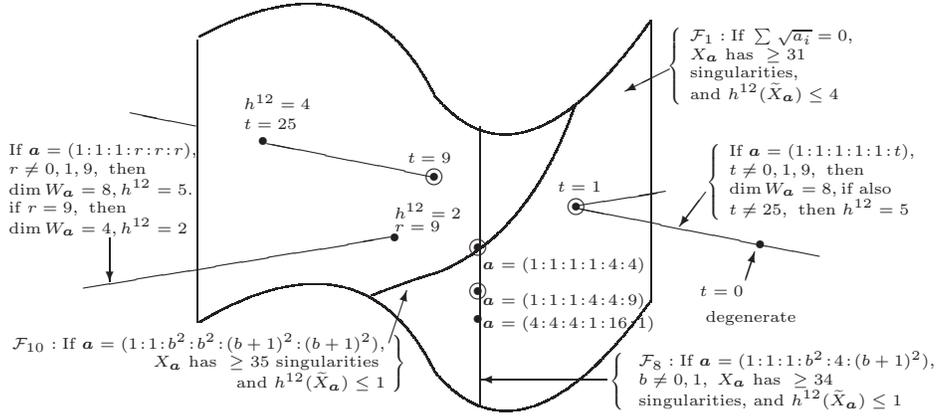
\begin{figure}[H]
$$
\setlength{\unitlength}{2173sp}%
%\begin{picture}(8847,4792)(-1042,-5441)
%\begin{picture}(8847,4792)(-1042,-5041)
\begin{picture}(8847,4492)(-1042,-5041)
\tiny
\put(4618,-3016){\circle*{100}}%
\put(4618,-3016){\circle{180}}%
{\put(1048,-2266){\circle*{100}}}%
{\put(2998,-2671){\circle*{100}}}%
{\put(2998,-2671){\circle{180}}}%
{\put(6718,-3451){\circle*{100}}}%
{\put(3500,-3480){\circle*{100}}}%
{\put(3500,-3480){\circle{180}}}%
{\put(3520,-2100){\line(0,-1){3200}}}
{\put(301,-1111){\line( 0,-1){3225}}}%
{\put(5476,-886){\line( 0,-1){3225}}}%
{\put(7393,-3586){\line(-5, 1){2763.461}}}%
{\put(2968,-2671){\line(-5, 1){1875}}}%
{\put(283,-2101){\line(-5, 1){750}}}%
{\put(-1000,-3950){\line(6,1){3500}}}
{\put(4640,-3010){\line(6,1){1000}}}
{\put(2550,-3367){\circle*{100}}}%
{\put(6553,-3856){\vector( 1, 3){106.500}}}%
{\put(5688,-1441){\vector(-2,-1){510}}}%
{\put(6100,-2846){\vector(-3,-4){300}}}%
{\put(2488,-4576){\vector( 1, 3){200}}}%
{
\qbezier(316,-1076)(2176,-76)(3000,-1750)
\qbezier(3000,-1750)(4000,-2950)(5460,-900)
}%
{\qbezier(316,-4336)(2176,-3236)(3000,-4900)
\qbezier(3000,-4900)(4000,-6100)(5460,-4100)
}%
{\qbezier(4600,-1866)(4068,-3371)(3003,-3781)
\qbezier(3003,-3781)(2893,-3796)(2280,-4031)
}%
\put(6098,-4336){degenerate}%
\put(5663,-1456){$\noteone$}%
\put(6108,-2776){$\notetwo$}%
\put(-2000,-2876){$\notefive$}%
\put(-700,-3376){\vector(0,-1){500}}
\put(6033,-4040){$t=0$}%
\put(-2000,-4851){$\notethree$}%
\put(4970,-4995){\vector(-1,0){1430}}
\put(5000,-5051){$\notefour$}%
\put(2700,-2521){$t=9$}
\put(833,-1906){$h^{12}=4$}%
\put(833,-2131){$t=25$}%
\put(2550,-3300){$r=9$}%
\put(2550,-3150){$h^{12}=2$}%
\put(4413,-2850){$t=1$}%
\put(3558,-3736){$\ba=(1\!:\!1\!:\!1\!:\!1\!:\!4\!:\!4)$}%
\put(3558,-4136){$\ba=(1\!:\!1\!:\!1\!:\!4\!:\!4\!:\!9)$}%
\put(3558,-4400){$\ba=(4\!:\!4\!:\!4\!:\!1\!:\!16\!:\!1)$}%
{\put(3500,-3980){\circle*{100}}}%
{\put(3500,-3980){\circle{180}}}%
{\put(3500,-4300){\circle*{100}}}%
\end{picture}
$$
\caption{Values of the parameter $\ba$ for
certain members and subfamilies of the family of Calabi-Yau threefolds $X_\ba$,
with $\ba$ for modular $X_\ba$ marked by a point, which is
circled if $X_\ba$ is rigid.
}
\label{subfamilies}
\end{figure}

Figure \ref{subfamilies} 
gives a schematic diagram of the $5$-dimensional family  $X_\ba$ which
we shall study, and some of its subfamilies
(a complete list is given in Table~\ref{subfamiliesTable}). 
The diagram gives the 
dimension of these subfamilies
and the value of $h^{12}$ of 
the big resolution $\widetilde X_\ba$ of the general member of the subfamily.  
Values of $\ba$ where $\widetilde X_\ba$ is
modular are marked with points, and those which are rigid with circled
points.  In \S~\ref{sec:elliptic_surfaces} we will see that for
two equal indices $a_i=a_j$, there is a corresponding elliptic surface
in $X_\ba$.  We call the piece of $H^3$ corresponding to these
elliptic surfaces $W_\ba$.  When $\dim H^{3}-\dim W_\ba=2,$ 
i.e., $2h^{12}=\dim W_\ba$, one
expects $X_\ba$ to be modular, which we will see
is the case for all $7$ marked points in the diagram.

Finally we would like to point out that the family of Calabi-Yau varieties which we are
considering in this paper has recently also appeared in a different context. C.~Borcea has 
studied these varieties in the context of configuration spaces of planar polygons (see \cite{Bo}
where these varieties are called {\em Darboux varieties}).

\subsection*{Acknowledgments}

We are grateful to the following institutions for support: to the DFG for
grant Hu 337/5-1 (Schwerpunktprogramm ``Globale Methoden in der
komplexen Geometrie'') and to the University of Essen 
and M.~Levine for hospitality
during a stay supported by a Wolfgang Paul stipend. We are also greatly
indebted to V.~Batyrev, W.~Fulton, J.~Koll\'ar and P.H.M.~Wilson 
whose comments on
intersection theory and Calabi-Yau manifolds were very helpful.
We also thank  N.~Fakhruddin 
for useful remarks on $\ell$-adic Galois representations and M.~Sch\"utt for
pointing out some misprints.

\subsection*{Notation}

In this paper we consider projective Calabi-Yau varieties defined by
polynomial equations with coefficients in $\ZZ$.  We work over the field
$k$, where $k=\CC$, $\QQ$, $\overline\QQ$, $\FF_p$ or $\overline\FF_p$.
Further notation is as follows.

\smallskip

\hspace{-0.8cm}
\begin{tabular}{ll}
$M_{A_4}$ & The $A_4$ root lattice, as a sublattice of $\ZZ^{5}$.\\
$N_{A_4}$ & $(M_{A_4})^\vee$, identified with a 
sublattice of $M_{A_4}\otimes\QQ$.\\
$\varepsilon_{ij}$ & Point in $M_{A_4}$ at $e_i-e_j$.\\
$\Delta_{A_4}$ & Polytope in $M_{A_4}\otimes\RR$ with vertices 
at $\varepsilon_{ij}$.\\%, i.e., the convex hull of the roots.\\
$\widetilde\Sigma_{A_4}$ & Fan in $N_{A_4}\otimes\RR$ given by all faces of
the Weyl chambers.\\
$\widetilde\Sigma^3$&$3$-dimensional cones in $\widetilde\Sigma_{A_4}$.\\
$\widetilde P$ & Smooth toric variety defined by $\widetilde\Sigma_{A_4}$.\\
$\AAA^4_\sigma$ & Affine piece of $\widetilde P$ corresponding to 
$\sigma\in\widetilde\Sigma^3$.\\
$T$ & Torus $(k^*)^4\subset \widetilde P$.\\
$T_\sigma$ &Orbit under the torus action corresponding to $\sigma\in
\widetilde\Sigma_{A_4}$;\\
& $\widetilde P=\bigsqcup_{\sigma\in\widetilde\Sigma} T_\sigma$.\\
$X_u$ & 
For $u\in\PP^{20}$, a 
Calabi-Yau threefold defined in $\widetilde P$
defined by $\Delta_{A_4}$.\\
$X_\ba$ & For $\ba\in\PP^{5}$, a
 member of a $5$ dimensional family of singular\\
& Calabi-Yaus in $\widetilde P$, with
$30$ nodes on $X_\ba\setminus T$
if $\prod_{i=1}^5 a_i\not=0$.\\
$\overline X_\ba$ & A Calabi-Yau given by
taking a
specific choice of small projective\\
& resolution of the $30$ singularities on $X_\ba\setminus T$.\\
$\widetilde{\overline X}_\ba$& Big resolution of remaining singularities on
$\overline X_\ba$.\\
$X_t$ & $X_{(1:1:1:1:1:t)}$.\\
$\widetilde X$ & Big resolution of $X=X_u, X_\ba,{\overline X}_\ba$
 or $X_t$, for $X$ irreducible.\\
$\widehat X$ & A choice of small projective 
resolution of $X$, if one  exists.
\end{tabular}

\section{Toric varieties}

%\subsection{A review of Batyrev's construction}

Batyrev \cite{Ba} constructs 
Calabi-Yau varieties as hypersurfaces in a toric variety defined by 
a reflexive polytope $\Delta$ (\cite[p.510]{Ba}), in a lattice $M$, as follows.

The pair $(M,\Delta)$ gives rise to a fan $\Sigma$ 
in the dual space $N_{\RR}=M_{\RR}^\vee$, 
and a strictly convex support
function $h$ on  $N_{\RR}$. 
Let $(P, {\cal O}_P(1))$ be the corresponding polarized toric
variety. In general $P$ is a singular Fano variety with Gorenstein
singularities and ${\cal O}_P(1)$ is the anticanonical bundle. 

Batyrev shows that the general
element in the anticanonical system is a Calabi-Yau variety, with canonical
singularities exactly at the singular points of $P$. A desingularization
$\widetilde P$ of $P$ can be constructed by taking the maximal projective
triangulation $\widetilde{\Delta^{\ast}}$ 
of the dual polytope $\Delta^{\ast}$. 
We denote the corresponding fan by $\widetilde{\Sigma}$.
The strict transforms of the anticanonical
divisors on
$P$ define a family of Calabi-Yau threefolds 
on $\widetilde P$ whose general element
$X_u$ is smooth.

\subsection{The polytopes $\Delta_{A_4}$, $\Delta_{A_4}^*$
and $\widetilde{\Delta_{A_4}^*}$}

We now describe the lattice and polytope to which we apply
Batyrev's construction.  
More generally, in the following construction
one may replace $A_4$ by $A_n$, defined similarly;  more details can be
found in \cite{DL},
\cite{procesi} and \cite{ludwig}.

Let $M$ be the $A_4$ root lattice, given  as the following % rank $4$ 
sublattice of $\ZZ^5$:
$$
M=M_{A_4}=
\left\{(x_1,x_2,x_3,x_4,x_5)\;\Big\vert\quad x_i\in\ZZ,\;\sum x_i=0\right\} 
\subset \ZZ^5.
$$
The inner product on $M$ is induced by the standard inner product on
$\ZZ^5$, and we identify the dual lattice $N:=(M_{A_4})^\vee$ with
a sublattice of $M_\RR:=M\otimes\RR$.

\begin{definition}
The polytope $\Delta=\Delta_{A_4}$ in $M_\RR$ 
is defined to be the convex hull of the
roots $\varepsilon_{ij}:=e_i-e_j, 1\le i, j\le 5,i\not=j$ of $M$,
where $e_1,\dots, e_5$ is the standard basis for $\ZZ^5$.
In \cite[p.427]{V} it is shown that $\Delta$ is reflexive.
\end{definition}

It is a simple combinatorial exercise to enumerate the faces of
$\Delta_{A_4}$, $\Delta_{A_4}^*$
and $\widetilde{\Delta_{A_4}^*}$.  We have the following results.

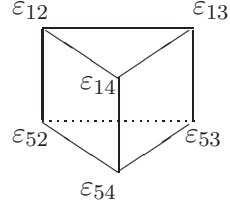
\begin{figure}
\begin{center} 
\unitlength 1mm 
\linethickness{0.15mm }
\begin{picture}(60,45)(0,0)
\thinlines
\put(0,5){\line( 5, 3){15}}
\put(0,5){\line(-2, 3){6}}
\put(0,5){\line(1, 6){3.52}}
\put(15,14.5){\line(-1, 1){11.5}}
\put(-6,14){\line( 3, 4){9}}
\qbezier[20](-6,14)(4.5,14.25)(15,14.5)
%\put(-6,14){\line( 5, 1){20}}
%          
\put(55,18.5){\line( -3,2) {9.5}}
\put(55,18.5){\line( 3, 2) {9.5}}
\put(65,25.1){\line(-1, 0){20}}
\put(55.2,6){\line( 0, 1){12.5}} 
\put(45,13){\line( 0,1){12}}
\put(65,12.8){\line( 0,1){12}}
\qbezier[20](45,12.8)(55,12.8)(65,12.8)
\put(55,6){\line(3,2) {9.5}}
\put(55,6){\line(-3,2){10}}
\put(-5,2){$\varepsilon_{14}$}
\put(-12,10){$\varepsilon_{15}$}
\put(-2,28){$\varepsilon_{12}$}
\put(12,10){$\varepsilon_{13}$}
\put(41,27){$\varepsilon_{12}$}
\put(41,10){$\varepsilon_{52}$}
\put(65,27){$\varepsilon_{13}$}
\put(64,10){$\varepsilon_{53}$}
\put(50,17){$\varepsilon_{14}$}
\put(50,3){$ \varepsilon_{54}$}
\put(-15,40){Tetrahedral face}
\put(-15,35){$F_1^+:=(x_1=1)\cap \Delta$}
\put(-15,-5){there are $10$ such faces}
\put(40,40){Prism face}
\put(40,35){$F_{15}^+:=(x_1 + x_5=1)\cap \Delta$}
\put(40,-5){there are $20$ such faces}
\end{picture} 
\end{center} 
\caption{Three dimensional faces of $\Delta$}
\label{deltas_three_faces}
\end{figure}

\begin{lemma}
\label{faces_of_delta_lem}
The polytope $\Delta$ has 
$20$ vertices, $60$ edges, 
$30$ square faces, $40$ triangular faces, and
$30$ three dimensional faces, given by 
$$
\begin{array}{lllllllll}
F_i^\varepsilon&:=& (x_i=\varepsilon1)\cap\Delta,&&
F_{ij}^\varepsilon &:=& (x_i+x_j=\varepsilon1)\cap\Delta,
\end{array}
$$
for $1\le i,j\le 5$, $i\not=j$,
and $\varepsilon=\pm$.
Two of these faces are shown in Figure~\ref{deltas_three_faces}.
\end{lemma}
\begin{proof}
See \cite[Lemma~1.18 and Korollar~1.19]{ludwig}.
\end{proof}

\begin{figure}
\setlength{\unitlength}{2368sp}%
\begin{picture}(8007,3559)(-2099,-9728)
\put(826,-9661){$=\frac{1}{5}(2,-3,2,2,-3)$}%
\put(-1574,-7486){$={F^+_{12}}^*$}%
\put(-1949,-9661){$=\frac{1}{5}(1,1,1,1,1,-4)$}%
\put(-2024,-7111){$\frac{1}{5}(3,3,-2,-2,-2)$}%
\put(5476,-7261){${F^+_{1}}^*$}%
\put(4351,-9661){${F^-_{5}}^*$}%
\put(-1274,-8761){${F^-_{35}}^*$}%
\put(1201,-7486){$={F^+_{1}}^*$}%
\put(-824,-9286){${F^-_{5}}^*$}%
\put(601,-7111){$-\frac{1}{5}(-4,1,1,1,1)$}%
\put(3676,-7111){${F^+_{12}}^*$}%
\put(  1,-8161){${F^-_{45}}^*$}%
\put(901,-8761){${F^+_{14}}^*$}%
\put(1501,-8161){${F^+_{13}}^*$}%
\put(1501,-9286){${F^-_{25}}^*$}%
\thinlines
\put(-194,-9331){\line( 0, 1){1425}}
\put(-719,-7261){\line( 4,-5){528}}
\put(877,-7261){\line( 4,-5){528}}
\put(1411,-9331){\line( 0, 1){1425}}
\put(-719,-8706){\line( 4,-5){528}}
\put(-734,-8701){\line( 0, 1){1425}}
\put(-179,-7921){\line( 1, 0){1575}}
\put(-179,-9346){\line( 1, 0){1575}}
\multiput(880,-8706)(117.33333,-146.66667){5}{\line( 4,-5){ 58.667}}
\multiput(871,-7246)(0.00000,-194.00000){8}{\line( 0,-1){ 97.000}}
\multiput(871,-8686)(-190.58824,0.00000){9}{\line(-1, 0){ 95.294}}
\put(-689,-7246){\line( 1, 0){1575}}
\put(5341,-7236){\line(-1,-2){1050}}
\put(4291,-9336){\line(-1, 4){525}}
\put(3766,-7236){\line( 1, 0){1575}}
\put(4291,-9321){\line( 0, 1){1425}}
\put(3766,-7251){\line( 4,-5){528}}
\put(5362,-7251){\line( 4,-5){528}}
\put(5896,-9321){\line( 0, 1){1425}}
\put(3766,-8696){\line( 4,-5){528}}
\put(3751,-8691){\line( 0, 1){1425}}
\put(4306,-7911){\line( 1, 0){1575}}
\put(4306,-9336){\line( 1, 0){1575}}
\multiput(5365,-8696)(117.33333,-146.66667){5}{\line( 4,-5){ 58.667}}
\multiput(5356,-7236)(0.00000,-194.00000){8}{\line( 0,-1){ 97.000}}
\multiput(5356,-8676)(-214.28571,0.00000){4}{\line(-1, 0){107.143}}
\multiput(3751,-8691)(250.00000,0.00000){2}{\line( 1, 0){125.000}}
\put(4306,-7896){\line( 3, 2){1020}}
\put(3811,-7296){\line( 1, 0){1395}}
\put(5176,-7341){\line(-1, 0){1335}}
\put(3886,-7401){\line( 1, 0){1155}}
\put(4936,-7476){\line(-1, 0){1005}}
\put(3991,-7551){\line( 1, 0){840}}
\put(4666,-7641){\line(-1, 0){615}}
\put(4141,-7731){\line( 1, 0){ 15}}
\put(4156,-7731){\line( 1, 0){ 30}}
\put(4186,-7731){\line( 1, 0){ 15}}
\put(4201,-7731){\line( 1, 0){ 15}}
\put(4216,-7731){\line( 1, 0){ 15}}
\put(4231,-7731){\line( 1, 0){ 15}}
\put(4246,-7731){\line( 1, 0){ 15}}
\put(4261,-7731){\line( 1, 0){ 15}}
\put(4276,-7731){\line( 1, 0){ 15}}
\put(4291,-7731){\line( 1, 0){ 15}}
\put(4306,-7731){\line( 1, 0){ 30}}
\put(4336,-7731){\line( 1, 0){ 15}}
\put(4351,-7731){\line( 1, 0){ 15}}
\put(4366,-7731){\line( 1, 0){ 15}}
\put(4381,-7731){\line( 1, 0){ 15}}
\put(4396,-7731){\line( 1, 0){ 15}}
\put(4411,-7731){\line( 1, 0){ 15}}
\put(4426,-7731){\line( 1, 0){ 15}}
\put(4441,-7731){\line( 1, 0){ 15}}
\put(4456,-7731){\line( 1, 0){ 15}}
\put(4471,-7731){\line( 1, 0){ 15}}
\put(4486,-7731){\line( 1, 0){ 15}}
\put(4501,-7731){\line( 1, 0){ 15}}
\put(4516,-7731){\line( 1, 0){ 15}}
\put(4531,-7731){\line( 1, 0){ 15}}
\put(4546,-7731){\line( 1, 0){ 15}}
\put(4561,-7731){\line( 1, 0){ 15}}
\put(4441,-7806){\line(-1, 0){225}}
\put(3871,-7371){\line( 0,-1){285}}
\put(3916,-7446){\line( 0,-1){435}}
\put(4006,-7551){\line( 0,-1){630}}
\put(4066,-7641){\line( 0,-1){765}}
\put(4141,-7731){\line( 0,-1){990}}
\put(4216,-7821){\line( 0,-1){1200}}
\put(3826,-7326){\line( 0,-1){165}}
\put(4306,-8121){\line( 1, 1){705}}
\put(4306,-8331){\line( 1, 1){960}}
\put(4291,-8541){\line( 1, 1){795}}
\put(4306,-8751){\line( 1, 1){540}}
\put(4291,-8961){\line( 1, 1){360}}
\put(4651,-8601){\line(-1,-1){210}}
\put(4291,-9111){\line( 1, 1){210}}
\put(5371,-8686){\line( 4,-5){108}}
\put(5533,-8911){\line( 4,-5){108}}
\put(5773,-9196){\line( 4,-5){108}}
\put(1288,-9196){\line( 4,-5){108}}
\put(1063,-8926){\line( 4,-5){108}}
\put(883,-8686){\line( 4,-5){108}}
\put(901,-8686){\line( 3,-4){495}}
\put(5356,-8671){\line( 4,-5){540}}
\put(-2099,-6361){Cubical face $\Theta_{15}=\varepsilon_{15}^*$ of $\Delta^*$}%
\put(3451,-6361){Face $\Theta_{12345}$ of $\widetilde{\Delta^*}$}
\put(4576,-8161){${F^-_{45}}^*$}
\end{picture}
\caption{Three dimensional faces of $\Delta^*$ and $\widetilde{\Delta^*}$}
\label{fig:faces_of_dual}
\end{figure}
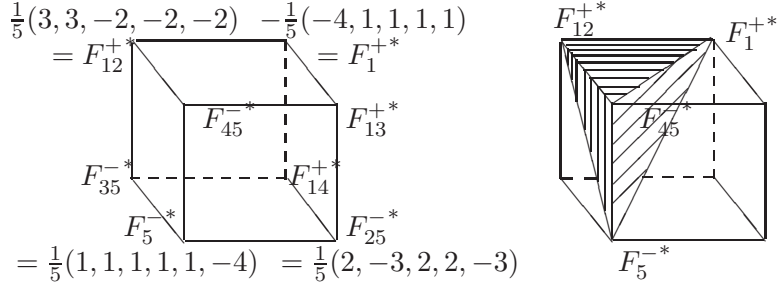

\begin{lemma}
\label{faces_of_delta_ast}
The dual polytope $\Delta^*$ has $30$ vertices and $20$ three dimensional
cubical faces, $\Theta_{ij}$ for $1\le i,j\le 5$,
$i\not=j$.  The vertices of $\Theta_{15}$ are shown in 
Figure~\ref{fig:faces_of_dual}, 
and  $\Theta_{ij}=\sigma\Theta_{15}$, where
$\sigma\in S_5$ with $\sigma:1,5\mapsto i,j$.
\end{lemma}
\begin{proof}
See \cite[Lemma~1.46]{ludwig}.
\end{proof}

\begin{lemma}
\label{faces_of_delta_ast_tilde}
The subdivded polytope $\widetilde{\Delta^*}$ has $120$ 
three dimensional faces, which are translations under $S_5$
of $\Theta_{12345}$ given in Figure~\ref{fig:faces_of_dual}.
\end{lemma}
\begin{proof}
See \cite[Beispiel~2.34]{ludwig}.
\end{proof}

In order to apply Batyrev's formulae for $h^{11}$ and $h^{12}$
we need to count the number of lattice points in the interior of the faces
of $\Delta$ and $\Delta^*$.  We make use of the following easy result.
\begin{lemma}
\label{points_in_cones}
If ${\bf w}_1,{\bf w}_2,{\bf w}_3,{\bf w}_4$ is a basis for
a lattice $L$, and $\Theta$ is a polytope with vertices
$\bf 0,{\bf w}_1,\dots,{\bf w}_4$, then
the only lattice points in $\Theta$ are its vertices.
\end{lemma}
By taking appropriate subdivisions of the faces $F^\varepsilon_{ij}$ and
$\Theta_{ij}$, we have
\begin{lemma}
\label{no_interior_points}
No proper face of $\Delta$ or $\Delta^*$
contains an interior lattice
point, and the only lattice point in the interior of 
$\Delta$ or $\Delta^*$ is the origin.
\end{lemma}

\subsection{The 
toric variety $\widetilde P$ defined by the fan $\widetilde\Sigma$}
\label{sub:toric_variety}

The fan
$\widetilde{\Sigma}$ in 
$N_{\RR}$ consists of cones given by the $120$ Weyl chambers
$$
\sigma_{ijklm}
=\left\{
(\alpha_1,
\alpha_2,
\alpha_3,
\alpha_4,
\alpha_5)\in\RR^5
\>\Big\vert\> \sum_{v=1}^5 \alpha_v = 0,\;
\alpha_i\ge \alpha_j\ge \alpha_k\ge \alpha_l\ge \alpha_m
\right\},
$$
where $\{i,j,k,l,m\}=\{1,2,3,4,5\}$, 
together with all their subfaces. 
E.g., $\sigma_{12345}$ is the cone on $\Theta_{12345}$ 
(see Figure~\ref{fig:faces_of_dual}).
The dual cones are given by
$$\sigma_{ijklm}^\vee:=\RR_{\ge0}(e_i-e_j) + 
\RR_{\ge0}(e_j-e_k) + 
\RR_{\ge0}(e_k-e_l) + 
\RR_{\ge0}(e_l-e_m).$$

We will consider Calabi-Yau threefold hypersurfaces in
the toric variety $\widetilde P$ defined by $\widetilde\Sigma$.
We first fix choices of
local and global coordinates for $\widetilde P$.

We identify the torus $T\cong (k^*)^4\subset \widetilde P$ with 
$\PP^4\setminus(\prod_{i=1}^5X_i=0)$, 
and use the  projective coordinates 
$X_1,\dots,X_5$ of $\PP^4$ when considering points in $T$.

For the affine piece
$\AAA^4_\varsigma:=\mathrm{Spec}(k[\varsigma(\sigma^\vee_{12345})])
\subset\widetilde P$, where $\varsigma\in S_5$,
we use coordinates $x_\varsigma,y_\varsigma,z_\varsigma, w_\varsigma$  
 corresponding to the basis $\varepsilon_{\varsigma(12)},
\varepsilon_{\varsigma(23)},\varepsilon_{\varsigma(34)},
\varepsilon_{\varsigma(45)}$ 
of $\varsigma(\sigma^\vee_{12345})$. Usually we just write $x,y,z,w$.

The identification of $T$
and $\AAA^4_\varsigma\setminus{(x_\varsigma y_\varsigma z_\varsigma w_\varsigma=0)}$ 
is given by
$$
\begin{array}{ccccccc}
x_\varsigma &=& X_{\varsigma(1)}/X_{\varsigma(2)},&&
y_\varsigma &=& X_{\varsigma(2)}/X_{\varsigma(3)},\\
z_\varsigma &=& X_{\varsigma(3)}/X_{\varsigma(4)},&&
w_\varsigma &=& X_{\varsigma(4)}/X_{\varsigma(5)}.
\end{array}
$$
This relationship between the coordinates of $\PP^4$ and of the
affine pieces of $\widetilde P$ is explained by the following lemma.
\begin{lemma}
\label{lem:cremona_transform}
The variety $\widetilde P$ is the graph of the Cremona transformation
$X_i\mapsto 1/X_i$ of $\PP^4$.
Thus $\widetilde P$ is obtained from $\PP^4$
by blowing up successively the (strict transforms of the) points 
$(1:0:0:0),(0:1:0:0)
\dots,(0:0:0:1)$, lines and planes spanned by any subset of
these points.  
\end{lemma}
\begin{proof}
See \cite[Lemma 5.1]{DL}.
\end{proof}

\subsubsection{Toric orbits in $\widetilde P$}
\label{description_of_P_minus_T}

There is a decomposition
$\widetilde P=\bigsqcup_{\sigma\in\widetilde\Sigma}T_\sigma$,
where $T_\sigma$ is the toric orbit of
$\widetilde P$ corresponding to $\sigma\in\widetilde\Sigma$.
Since $\widetilde \Sigma$ is given by taking cones on the faces of
$\widetilde\Delta^*$, we use the notation 
$T_\Theta:=T_{\RR_+\Theta}$
where $\Theta$ is a face of $\widetilde\Delta^*$.
By standard methods of toric geometry we have
\begin{eqnarray*}
\overline{T_{{F^\varepsilon_{i}}^*}}  &\cong&\widetilde P_{A_3},\\
\overline{T_{{F^\varepsilon_{ij}}^*}} &\cong&\widetilde P_{A_2}\times 
\widetilde P_{A_1}\cong
\widetilde\PP^2\times\PP^1,
\end{eqnarray*}
where $\widetilde P_{A_n}$ is the toric variety corresponding to the
root lattice $A_n$, and $\widetilde \PP^2$ is $\PP^2$ blown up in
$3$ points.
These are sketched in Figures~\ref{Sv_picture} and \ref{Pv_picture}.

\begin{figure}
\setlength{\unitlength}{2763sp}%
\begin{picture}(3637,3200)(1089,-2593)
\thinlines
{\put(4126,-436){\line(-1,-3){300}}
\put(3826,-1336){\line( 4,-5){300}}
\put(4126,-1711){\line( 1, 4){225}}
\put(4351,-811){\line(-3, 5){225}}
}%
{\put(1351,-1561){\line( 1, 1){750}}
\put(2101,-811){\line( 1,-1){750}}
\put(2851,-1561){\line(-1,-1){750}}
\put(2101,-2311){\line(-1, 1){750}}
}%
{\put(2101,-811){\line( 1, 3){300}}
}%
{
\put(3320,280){\line(-2, 1){300}}
\put(3000,440){\line(-5,-1){900}}
\put(2100,240){\line( 2,-1){300}}
\put(2400,85){\line( 5, 1){880}}
}%
{\put(4141,-436){\line(-6, 5){858}}
}%
{\put(3841,-1306){\line(-4,-1){992}}
}%
{\put(1101,-1171){\line( 2,-3){250}}
}%
{\put(1396,-466){\line( 1, 1){705}}
}%
{\put(2116,-2311){\line( 2,-1){540}}
\put(2656,-2581){\line( 3, 1){735}}
\put(3391,-2336){\line( 6, 5){750}}
}%
{\put(1396,-481){\line(-2,-5){270}}
}%
\put(2876,-661){$\widetilde\PP_2$}%
\put(2901,-2236){$\widetilde\PP_2$}%
\put(1651,-811){$\widetilde\PP_2$}%
\put(1576,-1636){$\PP^1\times\PP^1$}%
\put(2326,489){${}_{\PP^1\times\PP^1}$}%
\put(3926,-1261){${}_{\PP^1\times\PP^1}$}%
\put(4726, -1000)
{\parbox{2.4in}{\small 
There are $10$ copies of $\overline{T_{{F^\varepsilon_i}^*}}$ 
in $\widetilde P\setminus T$,
two being $\overline{\{x=0\}}$ and $\overline{\{w=0\}}$.
%The part at infinity of the toric 
The closures of $T_{\Theta^*}$, for $\Theta$ a two dimensional 
face of $F^\varepsilon_{j}$,
consists of $8$ copies of $\widetilde \PP^2$ and $6$ copies
of $\PP^1\times\PP^1$.  These intersect as indicated in this figure; 
hexagonal faces correspond to 
$\widetilde \PP^2$, and square faces correspond to $\PP^1\times\PP^1$.
}}
\end{picture}
\caption{The threefold $\overline{T_{{F^\varepsilon_i}^*}}$ 
in $\widetilde P\setminus T$.}
\label{Sv_picture}
\end{figure}

\begin{figure}
\setlength{\unitlength}{2763sp}%
\begingroup\makeatletter\ifx\SetFigFont\undefined%
\gdef\SetFigFont#1#2#3#4#5{%
  \reset@font\fontsize{#1}{#2pt}%
  \fontfamily{#3}\fontseries{#4}\fontshape{#5}%
  \selectfont}%
\fi\endgroup%
\begin{picture}(2475,2200)(401,-1123)
\thinlines
\put(1726,239){\line( 1, 0){975}}
\put(2701,239){\line( 6,-5){450}}
\put(3151,-136){\line(-2,-1){750}}
\put(2401,-511){\line(-1, 0){975}}
\put(1426,-511){\line(-6, 5){450}}
\put(976,-136){\line( 2, 1){750}}
\put(976,-136){\line( 0,-1){600}}
\put(1426,-511){\line( 0,-1){600}}
\put(2401,-511){\line( 0,-1){600}}
\put(3151,-136){\line( 0,-1){600}}
\put(3151,-736){\line(-2,-1){750}}
\put(2401,-1111){\line(-1, 0){975}}
\put(1426,-1111){\line(-6, 5){450}}
\put(1801,-211){$\widetilde\PP^2$}%
\put(1051,-661){${}_{\PP^1\times\PP^1}$}%
\put(1651,-811){${}_{\PP^1\times\PP^1}$}%
\put(2526,-661){${}_{\PP^1\times\PP^1}$}%
%\put(901,839){\parbox{11cm}{
\put(3576, -400){\parbox{7cm}{\small
\small 
There are $20$ copies of $\overline{T_{{F^\varepsilon_{ij}}^*}}$ 
in $\widetilde P$,
two being $\overline{\{y=0\}}$ and $\overline{\{z=0\}}$.
The closures of 
$T_{\Theta^*}$, for $\Theta$ a two dimensional 
face of $F^\varepsilon_{ij}$, consists of
$2$ copies of $\widetilde\PP^2$, and $6$ copies of $\PP^1\times\PP^1$,
corresponding to the hexagons and squares respectively in this figure.
}}%
\end{picture}
\caption{The threefold $\overline{T_{{F^\varepsilon_{ij}}^*}}$ 
in $\widetilde P\setminus T$.}
\label{Pv_picture}
\end{figure}

In terms of local coordinates $x,y,z,w$ for $\AAA^4_{\mathrm{id}}$, we have
$$
\begin{array}{lllll} 
{T_{{F^-_{4}}^*}} &=&
{T_{\frac{1}{5}(1,1,1,1,-4)}} &=  &{\{x=0\not=yzw\}},\\
{T_{{F^-_{45}}^*}} &=&
{T_{\frac{1}{5}(2,2,2,-3,-3)}}& = &{\{y=0\not=xzw\}},\\
{T_{{F^+_{12}}^*}} &=&             
{T_{\frac{1}{5}(3,3,-2,-2,-2)}}&= &{\{z=0\not=xyw\}},\\
{T_{{F^+_{1}}^*}} &=&              
{T_{\frac{1}{5}(4,-1,-1,-1,-1)}}&=&{\{w=0\not=xyz\}}.
\end{array}
$$
The intersections of the closures these hypersurfaces is sketched in 
Figure~\ref{oo_picture}.

\begin{figure}
\setlength{\unitlength}{2368sp}%
\begin{picture}(4296,3789)(4,-2400)
\thinlines
\put(958,818){\line(-6,-5){660}}
\put(298,268){\line(-1,-3){282}}
\put( 16,-578){\line( 4,-5){252}}
\put(283,-1492){\line( 1,-1){720}}
\put(1498,-1612){\line( 5,-1){615}}
\put(823,-937){\line( 0, 1){600}}
\put(1393,-472){\line( 1,-1){780}}
\put(1408,698){\line( 5, 1){570}}
\put(1978,812){\line( 5,-2){510}}
\put(3148,-1042){\line( 5,-6){290}}
\put(838,-352){\line( 1, 2){450}}
\put(3418,-112){\line(-4, 3){932}}
\put(2486,587){\line(-3,-1){687}}
\put(1799,358){\line(-3, 1){540}}
\put(1259,538){\line(-1, 1){301}}
\put(958,839){\line( 3,-1){420}}
\put(1378,699){\line( 6,-5){420}}
\put(1798,349){\line(-1,-2){411}}
\put(1387,-473){\line(-4, 1){520}}
\put(867,-343){\line(-1,-1){584}}
\put(283,-927){\line( 0,-1){565}}
\put(283,-1492){\line( 1, 1){555}}
\put(838,-937){\line( 1,-1){660}}
\put(1498,-1597){\line(-5,-6){500}}
\put(998,-2197){\line( 4,-1){516}}
%maybe following line overlaps too long line
\put(1514,-2326){\line( 1, 1){614}}
\put(2128,-1712){\line( 0, 1){490}}
\put(2128,-1222){\line( 5, 1){1005}}
\put(3133,-1021){\line( 1, 3){303}}
\put(3436,-112){\line( 3,-5){216}}
\put(3652,-472){\line(-1,-4){229}}
%following is too long:
\put(1860,-2200){\line(-5,3){280}}
\put(3423,-1388){\line(-6,-5){708}}
\put(2701,-1996){\line(-4,-1){840}}
\put(2140,-442){${}_{x=0}$}
\put(1030,-97){${}_{y=0}$}
\put(1105,-997){${}_{z=0}$}
\put(4300,-738){\parbox{6cm}{\small
The closures of $\{x=0\}$, $\{y=0\}$ and $\{z=0\}$ intersect as indicated
by the intersections of the corresponding polyhedra.
%The prisms correspond to
%$\{y=0\}$ and $\{z=0\}$, and the truncated octahrdra corresponds to 
%$\{x=0\}$.
A polyhedron corresponding to $\{w=0\}$ meets these polyhedra
in the labeled faces.  Where two polyhedra meet in a hexagon the corresponding
threefolds have intersection $\widetilde\PP^2$, and where they meet in a
square the corresponding threefolds intersect in $\PP^1\times\PP^1$.
}}
\end{picture}
\caption{How $\overline{\{x=0\}},
\overline{\{y=0\}}$ and
$\overline{\{z=0\}}\subset \widetilde P\setminus T$ meet.}
\label{oo_picture}
\end{figure}

\section{The Calabi-Yau varieties}

Following Batyrev, we define a 
family of hypersurfaces in $\widetilde P$,
given by elements $X\in|-K_{\widetilde P}|$.
The general member of
the family is given by an equation containing exactly 
the monomials corresponding to the lattice points of $\Delta$. 

In our case, $\Delta$ has $21$ lattice points, and
so gives a $20$ dimensional family.
The general member, when restricted to 
the open torus $T$, has an equation
\begin{equation}
\label{equation_for_whole_family}
X_u : \> \sum_{1\le i,j \le 5,i\not=j} u_{ij} X_iX_j^{-1}=t
\>\>\text{ for }u=(u_{12}:u_{13}:\dots:u_{45}:t)\in\PP^{20}
\end{equation}
in terms of the homogenous coordinates for $T\subset\PP^4$.

Given the above analysis of $\Delta$
and $\Delta^*$, we can now prove the following.

\begin{proposition}
\label{prop:values_of_hijXu}
For every smooth member $X_u$ of the family of Calabi-Yau threefolds
(\ref{equation_for_whole_family}), we have
\begin{enumerate}
\item[\rm{(i)}] The Euler number $e(X_u)=20$.
\item[\rm{(ii)}]The Hodge numbers of $X_u$ are given by
$h^{00}=h^{33}=h^{30}=h^{03}=1,\ h^{10}=h^{01}=h^{20}=h^{02}=0,\ h^{11}=26$ and
$h^{12}=h^{21}=16$.  
\end{enumerate}
\end{proposition}

\begin{proof}
Since $X_u$ is smooth and
Calabi-Yau the only Hodge numbers to be computed are $h^{11}$ and
$h^{21}$.  By \cite[p.521]{Ba} we have
$$
h^{11}(X_u) = l(\Delta^{\ast}) - 5 -
\sum_{\codim\Theta^{\ast}=1}l^{\ast}(\Theta^{\ast})
 + \sum_{\codim\Theta^{\ast}=2}l^{\ast}
(\Theta^{\ast})l^{\ast}(\Theta),
$$
$$
h^{21}(X_u) = l(\Delta) - 5 -
\sum_{\codim\Theta=1}l^{\ast}(\Theta)
 + \sum_{\codim\Theta=2}l^{\ast}
(\Theta)l^{\ast}(\Theta^{\ast}),
$$
where $l(\Theta)$ denotes the number of
lattice points in $\Theta$ and $l^{\ast}(\Theta)$
denotes the number of 
interior lattice points of $\Theta$,
for any face $\Theta$ of $\Delta$.

Lemmas~\ref{faces_of_delta_lem},
\ref{faces_of_delta_ast} and
\ref{no_interior_points} 
 imply that $l(\Delta)=21, l(\Delta^{\ast})=31$,
$l^*(\Delta^{\ast})=1$ and $l^*(\Theta^{\ast})=0$ for all proper faces 
$\Theta$ of $\Delta^{\ast}$. Hence
$
h^{11}(X_u)=26,\ h^{21}(X_u)=16.$
The Euler characteristic is then given by
$e(X_u)=2 h^{11}-2h^{12}=20.$
\end{proof}

\subsection{Resolution of singular Calabi-Yau threefolds}

We consider elements $X\in|-K_{\widetilde P}|$ which have $s$ nodes, but no other
singularities. We denote the big resolution of $X$,
obtained by blowing up the nodes, 
by $\widetilde X$. 
We also have $2^s$ small resolutions of $X$, where each node is replaced by a
$\PP^1$. 
It is not clear whether there are any small
projective resolutions as these could all contain null homologous lines.
By $\widehat X$ we denote a small
{\em projective} resolution, when one exists.

Let $X_u$ be a smooth member of the family 
(\ref{equation_for_whole_family}), and let $X_\ba$ be an element of the
family with $s$ nodes, but no other singularities.  Then we have

\begin{proposition} 
\label{h12_of_Xa}
Let $\widetilde h^{pq}=h^{pq}(\widetilde X_\ba)$, resp. $\widehat
h^{pq}=h^{pq}(\widehat X_\ba)$ be the Hodge numbers of the big resolution
$\widetilde X_\ba$ of $X_\ba$, resp. a small projective resolution of $X_\ba$. Then the following holds
\begin{enumerate}
\item [{\rm(i)}]
$ e(X_\ba)=e(X_u)+s, \quad e(\widehat X_\ba)=e(X_u)+2s,\quad
e(\widetilde X_\ba)=e(X_u)+4s$
\item [{\rm(ii)}]
$\widetilde h^{30}=\widetilde{h}^{03}=\hat h^{30}=\hat h^{03}=1$
\item [{\rm(iii)}]
$\tilde h^{10}=\tilde h^{01}=\hat h^{10}=\hat h^{01}=0, \quad
\tilde h^{20}=\tilde h^{02}=\hat h^{20}=\hat h^{02}=0$
\item [{\rm(iv)}]
$\tilde h^{11}-\tilde h^{12}=\frac 12 e(\widetilde X_\ba), 
\quad \hat h^{11}-\hat
h^{12}=\frac 12 e(\widehat X_\ba).$
\end{enumerate}
\end{proposition}

\begin{proof}
(i) These formulae are well known, cf. 
\cite[\S 1]{clemens} or \cite[Kapitel II]{We}.

(ii) Since $h^{p0}=h^{0p}$ are birational invariants, it
is enough to prove the assertion for $\tilde h^{30}$. Let $Q_1,\ldots, Q_s$ be
the exceptional quadrics 
in $\widetilde X_\ba$ and note that their normal bundle
is
$(-1, -1)$. Since $X_\ba$ is a Calabi-Yau 
variety with $s$ nodes it follows that
$\omega_{{\widetilde X}_\ba}={\cal O}_{\widetilde X_\ba}
\left(\sum\limits^s_{i=1} Q_i\right)$
and  $\tilde h^{30}=h^{0}(\omega_{\tilde X_\ba})=1$.

(iii) 
We consider the sequence 
$$0\rightarrow
\mathcal O_{\widetilde P}(-K_{\widetilde P})\rightarrow
\mathcal O_{\widetilde P}
\rightarrow
\mathcal O_{X_\ba}
\rightarrow
0.
$$
Since $h^1(\mathcal O_{\widetilde P})=0$ and 
$h^2(\mathcal O_{\widetilde P}(-K_{\widetilde P}))=
h^2(\mathcal O_{\widetilde P})=0$, we have
$h^1(\mathcal O_{X_\ba})=0$.  The resolution
$\pi:\widetilde X_\ba\rightarrow X_\ba$ is the blow up of double points,
i.e., of rational singularities, and hence
$R^1\pi_\ast\mathcal O_{\widetilde X_\ba}=0$.
By the Leray spectral sequence for the resolution this implies that
$h^1(\mathcal O_{\widetilde X_\ba})=h^1(\mathcal O_{X_\ba})=0$,
and hence, since these numbers are birational invariants,
$\tilde h^{10}=
\tilde h^{01}=
\hat h^{10}=
\hat h^{01}=0$.

It remains to prove that
$\tilde h^{02}=h^2({\cal O}_{\widetilde X_\ba})=h^1({\omega}_{\widetilde
X_\ba})=0$.
The latter can be deduced from the exact sequence
$$
0\rightarrow {\cal O}_{\widetilde X_\ba}\rightarrow {\omega}_{\widetilde X_\ba}
\rightarrow
\bigoplus\limits^s_{i=1} {\cal O}_{Q_i}(-1,-1) \rightarrow 0
$$
together with the fact that $H^1({\cal O}_{\widetilde X_\ba})=0$ which
we have already seen.

(iv) This follows immediately from (ii) and (iii).
\hfill
\end{proof}

\subsection{A five dimensional subfamily $X_\ba$ with $30$ nodes}
We now turn our attention to a certain subfamily of
(\ref{equation_for_whole_family})
of Calabi-Yau threefolds of the form
\begin{eqnarray}
\label{big_subfamily}
&X_{\boldsymbol{a}}\cap T:&
(X_1+\ldots +X_5)\left(\frac {a_1}{X_1}+\ldots +\frac {a_5}{X_5}\right)=a_6=t\\
&&\hspace{2cm}
\>\>\text{ for } \boldsymbol a = (a_1:a_2:a_3:a_4:a_5:t)\in\PP^5.
\nonumber
\end{eqnarray}
The variety $X_\ba$ is the closure of
$X_\ba\cap T$ in the toric variety $\widetilde P$. 
For $t\in\CC$, we will also use the notation $X_t:=X_{(1:1:1:1:1:t)}.$

In terms of the local coordinates
$x=X_1/X_2, y=X_2/X_3, z=X_3/X_4, w=X_4/X_5$
for $\AAA^4_{id}$, given in \S~\ref{sub:toric_variety},
the equation (\ref{big_subfamily}) for $X_{\boldsymbol a}$ becomes
\begin{equation}
\label{big_local_equation}
(a_1 + a_2x + a_3xy + a_4xyz + a_5xyzw)(1 + w + wz + wzy + wzyx) = a_6xyzw.
\end{equation}
\begin{remark}
If $\ba=(a_1:\dots :a_6)\in\PP^5$ with $\prod_{i=1}^6a_i=0$,
then $X_\ba$ is not irreducible, and so in general we assume
$\prod_{i=1}^6a_i\not=0$.
\end{remark}
Other local equations are given by permuting the $a_i$
appropriately.  Although this equation is not symmetric 
in the $a_i$, we will now see that up to birational equivalence
$a_6$ plays the same role as the other $a_i$.

\begin{lemma}
The variety $X_\ba$ defined by (\ref{big_subfamily}) 
is birational to a variety 
in $\PP^5$ defined by two equations:
\begin{equation}
\label{def_by_two_equations}
\sum_{i=1}^6 \frac{a_i}{X_i}=\sum_{i=1}^6 {X_i} = 0.
\end{equation}
\end{lemma}
\begin{proof}
This follows immediately from setting $X_6=-\sum_{i=1}^5 X_i$.
\end{proof}
\begin{corollary}
For any permutation $\sigma\in S_6$, the
varieties $X_\ba$ and $X_{\sigma(\ba)}$ are birational.
\end{corollary}

\begin{remark}
The Barth-Nieto quintic $N_5$ is the variety in $\PP^5$ defined by
(\ref{def_by_two_equations}) with all $a_i=1$.  A corollary of
the above lemma is that
$$X_{(1:1:1:1:1:1)}=X_1\sim_{\mathrm{bir}} N_5.$$
\end{remark}

\subsubsection{The singularities of $X_\ba$ on $X_\ba\cap T$}

\begin{lemma}
\label{sings_on_T}
For $\boldsymbol{a}=(a_1:a_2:a_3:a_4:a_5:t)\in\PP^5$, with $t\not=0$,
the variety $X_{\boldsymbol a}$ (over any field)
has a singularity at ${\boldsymbol b}\in T$ if and only if
$\boldsymbol a = \phi(\boldsymbol b)$ for some $\boldsymbol b\in T\subset
\PP^4$,
where $\phi$ is the map
\begin{eqnarray}
\label{definition_of_phi}
\phi : T &\rightarrow & \PP^5\\
(a:b:c:d:e) &\mapsto & (a^2:b^2:c^2:d^2:e^2:(a+b+c+d+e)^2).
\nonumber
\end{eqnarray}
\end{lemma}
\begin{proof}
Writing (\ref{big_subfamily}) 
as $f\cdot g=t$, and differentiating with respect to
$X_i$ gives
\begin{eqnarray}\label{eqn2}
g-\frac {a_i}{X^2_i} f=0,\quad i=1,\ldots, 5,
\end{eqnarray}
which implies that  a singular point has the form
$P=(\pm \sqrt{a_1}:\ldots: \pm \sqrt{a_5})$. 
Substituting  into
(\ref{big_subfamily}), we obtain
$$(\sqrt{a_1} \pm \sqrt{a_2} \pm \sqrt{a_3} \pm \sqrt{a_4} 
\pm \sqrt{a_5})^2 = t,$$
and so $\ba$ has the claimed form.
\end{proof}

\begin{proposition}
\label{A1_sings}
The singularity $\bb=(b_1:b_2:b_3:b_4:b_5)\in X_{\phi(\bb)}\cap T
\subset\PP^4$ 
in Lemma~\ref{sings_on_T} is an $A_1$ singularity.
\end{proposition}
\begin{proof}
Consider the affine cover of $T$
with coordinates $y_i = X_i/X_5-b_i/b_5$ for $1\le i\le 4$.  
In terms of these coordinates, the singularity is at $(0,0,0,0)$, and
the equation for $X_{\phi(\bb)}$ is given by
\begin{eqnarray}
\label{equation_near_singularity}
(b_1+b_2 + b_3 + b_4 + b_5)
\left(\frac{y_1^2}{b_1}
+\frac{y_2^2}{b_2}
+\frac{y_3^2}{b_3}
+\frac{y_4^2}{b_4}
\right)&&\nonumber
\\
- (y_1 + y_2 + y_3 + y_4)^2
+\text{higher order terms} &=& 0.
\end{eqnarray}
The matrix of the quadratic form given by the degree two part of
(\ref{equation_near_singularity}) is 
\begin{equation}
\label{quadratic_form}
M:=\left(
\begin{array}{cccc}
s/b_1 -1 & -1 & -1 & -1 \\
-1 & s/b_2-1 & -1 & -1\\
-1 & -1 &s/b_3-1 & -1\\
-1 & -1 & -1 &s/b_4-1
\end{array}
\right),
\end{equation}
where $s=b_1 + b_2 + b_3 + b_4 + b_5$.  We have
$$\det(M)=\frac{b_5(b_1+b_2 + b_3 + b_4 + b_5)^3}{b_1b_2b_3b_4},$$
which is non zero since $\bb\in T$.  Thus the singularity is as claimed.
\end{proof}

\subsubsection{The singularities of $X_\ba$ on $X_\ba\setminus T$}

\begin{lemma}
\label{local_sing_not_on_T}
Let $\ba=(a_1:\cdots:a_6)\in\PP$ with all $a_i\not=0$, and $\varsigma\in S_5$.
Then over any field,
$X_{\boldsymbol a}$ has singularities in 
$(X_{\boldsymbol a}\setminus T) \cap {\mathbb A}^4_\varsigma$ only at the point
$$(x_\varsigma,y_\varsigma,z_\varsigma,w_\varsigma)=
(-a_{\varsigma(1)}/a_{\varsigma(2)},0,0,-1),$$
where $x_\varsigma,y_\varsigma,z_\varsigma,w_\varsigma$ are as in
\S~\ref{sub:toric_variety}.
This singularity is a node.
\end{lemma}
\begin{proof}
This follows easily by computating the partial derivatives of
(\ref{big_local_equation}), and setting them to zero.
\end{proof}

\begin{lemma}
\label{exactly_30_sings}
For general $\ba\in\PP^5$ there are
exactly
 $30$ singularities on $X_{\boldsymbol a}$.
\end{lemma}
\begin{proof}
By Lemma~\ref{sings_on_T}, in general, $X_\ba\cap T$ is smooth.
By Lemma~\ref{local_sing_not_on_T}, there is only one singularity in 
the affine piece ${\mathbb A}^4_\varsigma$,
contained in 
$
(y_\varsigma=z_\varsigma=0\not=x_\varsigma w_\varsigma)=
T_{\varsigma(F^+_{12}\cap F^-_{45})^*}
$.
Since $\Delta$ has
${5\choose2}{3\choose2}=30$ faces 
$\sigma(F^+_{12}\cap F^-_{45})$ for $\sigma\in S_5$,
this implies the result.
\end{proof}

\subsubsection{Classification of singular subfamilies}

By Lemma~\ref{sings_on_T}, to find the number 
of singularities on $X_\ba\cap T$, we need
to determine the number of $\bb$ with $\ba=\phi(\bb)$.  We obtain
the following result.

\begin{proposition}
\label{prop:number_of_nodes_on_subfamilies}
For $\bb=(b_1:\cdots:b_5)\in T$,
the number of nodes on $X_\phi(\bb)$ is given by
\begin{equation}
\label{eqn:more_extra_sings}
30 + 
\#\left\{J : J\subset \{1,2,3,4,5\},\; \sum_{i\in J}b_i=0\right\}.
\end{equation}
\end{proposition}
\begin{proof}
Lemma~\ref{exactly_30_sings} gives $30$ nodes on $X_\ba\setminus T$.  
By Lemma~\ref{sings_on_T}, the only nodes on $X_\ba\cap T$ are at
$\bc$ with $\phi(\bc)=\phi(\bb)$.
Let $\bc=(c_1:\cdots:c_5)\in\PP^4$.
After scaling, $b_i=\eps_i c_i$, with $\eps_i=\pm 1$
and $\sum b_i = \sum c_i$.  From this we obtain
$$\sum (1-\eps_i)b_i=0.$$
Thus the correpondence between nodes in $T$ 
and subsets of $\{1,2,3,4,5\}$ is given by
\begin{equation}
\label{eqn:sign_changes}
\bc \longleftrightarrow \{i : c_i\not=b_i\},
\end{equation}
with the empty set corresponding to $\bb$.
\end{proof}

\begin{remark}
\label{rem:notes_on_table}
By considering how subsets of the $b_i$ can intersect,
this lemma allows us to classify all subfamilies of $X_\ba$ with more
than $30$ nodes.
These are given in Table~\ref{subfamiliesTable}.
In this table, sets $\mathcal F_i$ are given in a shortened form, e.g., 
$\mathcal F_7=\phi\{(a:a:a:-a:b) \}$ should be read as
$$
\mathcal F_7 := \phi\{(a:a:a:-a:b)\in\PP^4\; | \; a,b\in k\setminus \{0\}\}.
$$
To use this table, one must take the smallest set $\mathcal F_i$ containing
a given $\ba$, up to permutation of coordinates; e.g.;
$(1:1:1:1:1:1)\in 
\mathcal F_{15}\subset
\mathcal F_{11}\subset\mathcal F_{4}\subset\mathcal F_{2}\subset\mathcal F_{0}$; the
data for $X_1$ is given by the last line of the table.
\end{remark}

\begin{table}%[H]
$$
\begin{array}{lrlclll}
\text{Smallest family } \mathcal F_i&\multicolumn{2}{r}{\text{dimension}} & 
\multicolumn{1}{l}{\text{number of}}
& e(\widetilde{ X}_\ba)
& e(\widetilde{\overline X}_\ba) & e(\widehat{ X}_\ba)\\
\text{containing } \ba  &\multicolumn{2}{r}{\text{of }\mathcal F_i}
& \text{nodes on } X_{\ba}  
&&&\\
\multicolumn{2}{l}{\mathcal F_0=\{\ba\in\PP^6, a_i\not=0 \}     }& 5 & 30 & 140 & 80  & 80 \\   
\multicolumn{2}{l}{
\mathcal F_1=\phi\{(a:b:c:d:e)\}   }       & 4 & 30+1 & 144 &  84  & -\\
\multicolumn{2}{l}{\mathcal F_2=\phi\{(a:-a:b:c:d)\}}      & 3 & 30+2 &148&88& -    \\   
\multicolumn{2}{l}{\mathcal F_3=\phi\{(a:b:-a-b:c:d)\}}    & 3 & 30+2 &148&88& -    \\
\multicolumn{2}{l}{\mathcal F_4=\phi\{(a:-a:a:b :c)\}}     & 2 & 30+3 &152&92& -     \\
\multicolumn{2}{l}{\mathcal F_5=\phi\{(a:-a:b:a-b:c)\}}    & 2 & 30+3 &152&92& -     \\
\multicolumn{2}{l}{\mathcal F_6=\phi\{(a:-a:b:-b:c)\}}     & 2 & 30+4 &156&96&88  \\
\multicolumn{2}{l}{\mathcal F_7=\phi\{(a:a:a:-a:b) \}  }   & 1 & 30+4 &156&96& -     \\   
\multicolumn{2}{l}{\mathcal F_8=\phi\{(a:a:a:-2a:b)\} }    & 1 & 30+4 &156&96& -     \\   
\multicolumn{2}{l}{\mathcal F_9=\phi\{(a:a:b:-b:b-a)\}}    & 1 & 30+4 & 156&96& -   \\    
\multicolumn{2}{l}{\mathcal F_{10}=\phi\{(a:a:b:b:-a-b)\}}    & 1 & 30+5 & 160 &  100  & - \\        
\multicolumn{2}{l}{\mathcal F_{11}=\phi\{(a:a:-a:-a:b)\}}     & 1 & 30+6 & 164 &  104  & 92   \\
\multicolumn{2}{l}{\mathcal F_{12}=\phi\{ (1:1:1:1:-1)\}} &0  &30+5 & 160 &  100  & - \\ 
\multicolumn{2}{l}{\mathcal F_{13}=\phi\{ (1:1:1:2:-2)\}} &0 &30+5 & 160 &  100  & - \\
\multicolumn{2}{l}{\mathcal F_{14}=\phi\{ (1:1:1:1:-2)\}} &0  &30+7 & 168 &  108  & - \\
\multicolumn{2}{l}{\mathcal F_{15}=\phi\{ (1:1:1:-1:-1)\}} &0 &\phantom{2}30+10& 180 &  120 & 100
\end{array}
$$
\caption{
Number of nodes on $X_\ba$ for 
subfamilies of $\ba$, up to permutation of coordinates,
as images of subfamilies of $\bb\in T$, under $\phi$, given by
(\ref{definition_of_phi})
; see Remark~\ref{rem:notes_on_table}.  The last $4$ sets are
$\{(1:1:1:1:1:9)\}, \{(1:1:1:4:4:9)\}, \{(1:1:1:1:4:4)\}$ and 
$\{(1:1:1:1:1:1)\}$ respectively.
}
\label{subfamiliesTable}
\end{table}

As a corollary of Proposition~\ref{prop:number_of_nodes_on_subfamilies},
 we have the following result.

\begin{proposition}
\label{lem:euler_numbers}
The Euler numbers of $\widetilde X_{\ba}$,
$\widetilde{\overline X}_{\ba}$, and of the small resolution
 $\widehat{ X}_{\ba}$, if it exists, are given in Table~\ref{subfamiliesTable}.
\end{proposition}
\begin{proof}
This follows 
from Proposition~\ref{prop:values_of_hijXu}(i) 
and  Proposition~\ref{h12_of_Xa}(i).
\end{proof}

We will discuss whether $\widehat{X}_{\ba}$ exists
 in \S~\ref{sub:small_resolution}, \S~\ref{subsec:smallres_in_T}
and \S~\ref{subsec:existance_of_small_res}.

\begin{remark}
In Proposition~\ref{prop:values_of_h12} we will see that $h^{12}$
of the general member of any of the families in 
Table~\ref{subfamiliesTable} is equal to the
dimension of the family.  In particular, we will see
that for $\ba$ in a zero dimensional family, $X_\ba$ is rigid.
\end{remark}

\subsubsection{Weil divisors and small resolutions}
\label{sub:small_resolution}

We shall now discuss whether $X_\ba$ has a small projective
resolution. If this is the case then $X_\ba$ possesses a smooth
projective Calabi-Yau model. We shall see later 
(Proposition~\ref{existence_of_smooth_models}) that
these two statements are in fact equivalent.

First we recall a fact about small resolutions.
Let $X$ be a projective variety, and $P\in X$ a  double point.
Near $P$, 
 $X$ is locally analytically
a cone over a quadric surface $Q$.
The local divisor class group is generated by the cones 
$\Lambda_1, \Lambda_2$ over the rulings of $Q$.
Blowing up $X$ in $P$ defines a (big) 
resolution $\pi:\widetilde X\rightarrow X$, with
$\pi^{-1}(P)\cong Q$.
 Locally analytically
we can also blow up $\Lambda_1$ and $\Lambda_2$, obtaining
small resolutions $p_i:X_i\rightarrow X$, $i=1,2$, with
$p_i^{-1}(P)\cong\PP^1$.  However, this is an analytic construction, and it
is not clear whether the $X_i$ are projective.  We have a 
commutative diagram
$$
\xymatrix{
&\widetilde X\ar[dr]^{q_2}\ar[dl]_{q_1}\ar[dd]_\pi\\
X_1\ar[dr]_{p_1}&&X_2\ar[dl]^{p_2}\\
&X
}
$$
where the $q_i$ each blow down one family of rulings.
This is Atiyah's flop.
If $W$ is a global Weil divisor on $X$ through $P$,
blowing up along $W$ gives a projective variety.
If $W$ is not Cartier near $P$
then it is  analytically locally equivalent to $a\Lambda_1$ or
$a\Lambda_2$ for some $a>0$.
By the universal property of the blow up, blowing up $X$ along
$W$ is the same as performing a small resolution.
For a discussion of this material see 
\cite[example IV-27]{EH}.

\subsubsection{Small resolution 
 $\overline X_{\ba}$ of the $30$ singularities on
$X_{\ba}\setminus T$}

We now show that 
we can find a projective variety which is a small projective resolution of
each of the $30$ nodes on $X_\ba\setminus T$.
First we will define surfaces which define the Weil divisors we use
for the blow up, as described above.

\begin{definition}
\label{def:surfacesSij}
Let $1\le i<j\le 5$, $i\not=j$ and $\varepsilon=\pm$.
For some $\varsigma\in S_5$, we have
$\overline{T_{{F^\varepsilon_{ij}}^*}}
=\overline{(y_\varsigma=0)}$ in the affine piece $\AAA^4_\varsigma$.
We define the surface $S^{\varepsilon ij}_\ba$ 
to be one of the two components of 
$\overline{(y_\varsigma=0)}\cap X_\ba$,
given in terms of the coordinates for $\AAA^4_\varsigma$ by
\begin{equation}
\label{eq:set_y=0}
S^{\varepsilon ij}_\ba : \;
\overline{(a_{\varsigma(1)} + a_{\varsigma(2)}x_\varsigma= y_\varsigma=0)}.
\end{equation}
This is independent of which of the $6$ possible $\varsigma$ is chosen.
The surface $S^{-45}_{\ba}$ is indicated in Figure~\ref{intersect_picture}.
\end{definition}

\begin{lemma}
The $10$ surfaces $S_{\ba}^{-ij}$ are smooth, disjoint, and
each contains $3$ nodes of $X_\ba$.
\end{lemma}
\begin{proof}
Smoothness is clear.
Set $T_{ij}:=T_{{F^-_{ij}}^*}$.
For $\{i,j\}\not=\{k,l\}$
the vertices ${F^-_{ij}}^*$ and ${F^-_{kl}}^*$  do not 
lie on a common edge of $\Delta^*$ (see Figure~\ref{fig:faces_of_dual}),
so
$\overline{T_{ij}}\cap \overline{T_{kl}} =\varnothing$,
from which disjointness follows.
From the defining equation, $S_{\ba}^{-ij}$ 
contains the node $(-a_{\varsigma(1)}/a_{\varsigma(2)},0,0,-1)$
in $(X_\ba\setminus T)\cap\AAA^4_\varsigma$.
As in Lemma~\ref{exactly_30_sings}, each node lies on a surface
$T_{(F^+_{ij}\cap F^-_{kl})^*}$.  Since
$T_{ij}$ contains $3$ such surfaces, $S_{\ba}^{-ij}$ contains $3$ nodes.
This can be seen in Figure~\ref{intersect_picture}.
\end{proof}

\begin{remark}
Since $S_{\ba}^{\varepsilon ij}\cap T
\subset\overline{T_{{F^\varepsilon_{ij}}^*}}\cap T=\varnothing$, 
the surfaces $S_{\ba}^{\varepsilon ij}$ do not pass through any 
singularities $\bb\in X_\ba\cap T$, 
and may be ignored in considering
resolutions of such singularities.
\end{remark}

\begin{definition}
We define $\overline{X_\ba}$ to be the blow up of $X_\ba$ along
all $10$ surfaces $S^{-ij}_\ba$.
This is a  projective small resolution of all
 $30$ nodes in $X_\ba\setminus T$.
\end{definition}

\begin{proposition}
\label{prop:smooth_family}
Let $\ba\neq\phi(\bb)$ for any $\bb\in\PP^4$. 
Then the nodal Calabi-Yau variety $X_\ba$ has a
small projective resolution 
$\overline X_{\boldsymbol a}$.
\end{proposition}
\begin{proof}
This is because $X_\ba\cap T$ is smooth in this case,
by Lemma~\ref{sings_on_T}.
\end{proof}

We have now constructed a five dimensional family of
Calabi-Yau threefolds,  $\overline X_{\boldsymbol a}$.
In Proposition~\ref{prop:values_of_h12} we will see
that $h^{12}$ of the general member is $5$.
By Proposition~\ref{h12_of_Xa} 
the Euler number of $X_\ba$ for general $\ba$ is
$80$, and so the Hodge diamond is as follows:
$$
\begin{array}{ccccccc}
 & &&1&&&\\
 & &0&&0&&\\
 &0&&45 & &0&\\
1& &5&&5&&1\\
 &0&&45&&0&\\
 & &0&&0&&\\
 & &&1&&&\\
\end{array}
$$

\subsubsection{Big resolutions of singularities on $X_\ba\cap T$}

As discussed in \S~\ref{sub:small_resolution},
the big resolution of a singularity $\bb\in X_\ba\cap T$
replaces $\bb$ by a quadric, which in this case has
equation  given by the quadratic
part of (\ref{equation_near_singularity}).
By transforming (\ref{quadratic_form}) to a diagonal matrix,
one can see that the rulings of this quadric are defined over the field
$$\QQ\left[\sqrt{-b_0b_1b_2(b_0 + b_1 + b_2)},
\sqrt{-b_3b_4b_5(b_3+b_4+b_5)}\right],$$
and other similarly defined fields.
For later use, we now consider the number of points on these quadrics
over finite fields.
First we need the following:
\begin{lemma}
Let $f$ be an irreducible quadric in $\PP^3(\FF_p)$, for $p\not=2$,
with corresponding matrix $M\in M_4(\FF_p)$.  Then we have
$$\#(f=0)(\FF_p)=\left\{
\begin{array}{ll}
(p+1)^2 & \text{ if }\det{M}\text{ is a square in } \FF_p \\
p^2 + 1 & \text{ if }\det{M}\text{ is not a square in } \FF_p. 
\end{array}
\right.
$$
\end{lemma}
\begin{proof}
This follows from  
\cite[Proposition~5, IV\S1.7]{serre}.
\end{proof}
Together with
 Proposition~\ref{A1_sings}, this immediately implies the following:
\begin{corollary}
\label{cor:nodes_of_X}
For the quadric $Q_\bb$ introduced in resolving the singularity at 
$\bb:=(b_1:b_2:b_3:b_4:b_5)$ on $X_{\phi(\bb)}$, we have
$$
\#Q_\bb(\FF_p)
=\left\{
\begin{array}{ll}
(p+1)^2 & \text{ if }\sqrt{b_0b_1b_2b_3b_4b_5}\in\FF_p \\
p^2 + 1 & \text{ if }\sqrt{b_0b_1b_2b_3b_4b_5}\notin\FF_p,
\end{array}
\right.$$
where $b_0=b_1+b_2+b_3+b_4+b_5$.
\end{corollary}

\subsubsection{Small resolutions of singularities on $X_\ba\cap T$}
\label{subsec:smallres_in_T}

\begin{proposition}
\label{prop:case_of_small_res}
For $\bb=(a : b : -b : c : -c)$ (up to permutation of coordinates), there
is a small projective resolution of $\bb\in X_{\phi(\bb)}$.
\end{proposition}
\begin{proof}
In this case 
the surface $(X_2+X_3=X_4+X_5=0)$ lies in $X_{\phi(\bb)}$,
contains the point $\bb$, and its closure is a smooth Weil divisor in
$X_{\phi(\bb)}$.  Thus a small resolution is achieved by
blowing up this surface.
\end{proof}

\begin{remark}
Note that for $\bb$ as in the proposition, $X_{\phi(\bb)}$ has at 
least four singularities, at
$(a : b : -b : c : -c), (a : -b : b : c : -c),
(a : -b : b : -c : c)$ and $(a : b : -b : -c : c)$.
\end{remark}

\begin{corollary}
\label{cor:singsofX1}
There is a small projective resolution of
all the singularities of $X_1$.
\end{corollary}
\begin{proof}
This follows from Proposition~\ref{prop:case_of_small_res}, since
if $\ba=(1:1:1:1:1:1)=\phi(\bb)$ then $\bb=(1:1:1:-1:-1)$,
up to permutation of coordinates.
\end{proof}

\begin{remark}
\label{rem:explicit_small_res}
The only Calabi-Yau threefolds in Table~\ref{subfamiliesTable}
for which all singularities
can be resolved by the small resolution of
Proposition~\ref{prop:case_of_small_res}
are $X_{\phi(\bb)}$ for $\bb=(1:-1:b:-b:c)$ up to permutation,
with $1\pm b\pm c\not=0$ for all sign choices; this can easily be seen by
listing all possible singularities, by using (\ref{eqn:more_extra_sings}) and
(\ref{eqn:sign_changes}).
In Corollary~\ref{cor:examples_with_no_smoothCY}
we will see that other $X_\ba$ have no small projective resolution.

\end{remark}

\section{$X_\ba$ as a fibre product}
\label{sec:Xaasafibreproduct}

We now show that $X_\ba$ is
birational to the fibre product of families of elliptic curves.
This enables us to apply Schoen's results \cite{schoen} on such threefolds.

First we define the families of elliptic curves.

\begin{definition}
\label{def:elliptic_curve_family}
For $a,b,c\in k\setminus\{0\}$ define an elliptic 
surface $\mathcal E_{a,b,c}$ to be the resolution
of the surface $\mathcal E_{a,b,c}'\subset\PP^2\times\PP^1$
given by
\begin{eqnarray}
\label{eqn:Elliptic_abct_family}
\mathcal E_{a,b,c}':&&(x+y+z)(axy + byz + czx)t_0=t_1xyz,
\end{eqnarray}
where $(x:y:z)\in\PP^2, (t_0:t_1)\in\PP^1$.
There is a projection
$p:\mathcal E_{a,b,c}\rightarrow\PP^1$, with
fibre $\mathcal E_{a,b,c,t}:=p^{-1}(1:t)$.
We write $\mathcal E_{a}:=\mathcal E_{1,1,a}$
and $\mathcal E_{a,t}:=\mathcal E_{1,1,a,t}$.
\end{definition}

The only singular points of 
$\mathcal E_{a,b,c}'$ are the three singular points of
the fibre at infinity.  When these are resolved
the fibre at infinity becomes an $I_6$ fibre.

When $\mathcal E_{a,b,c,t}$ is smooth,
taking the zero to be $(0:1:-1)$, the elliptic curve
$\mathcal E_{a,b,c,t}$ is isomorphic to the cubic curve with equation
\begin{equation}
\label{eqn:weierstrassform}
\!y^2 = x\!
\left(\!\!\!
\left(\!
x + \frac{2(t^2 + a^2 + b^2 + c^2) \!-\! (t+a+b+c)^2}{8a^2}\!\right)^2 
-
\frac{A(a,b,c,t)}{64a^4}
\!\!\!\right),\!\!
\end{equation}
where
$$A(a,b,c,t):=
\prod_{(\nu,\mu)\in\{-1,1\}^2}
\left(t-(\sqrt{a}+\nu\sqrt{b}+\mu\sqrt{c})^2\right).$$
From this we find that the $j$-invariant of $\mathcal E_{a,b,c,t}$ is given by
\begin{equation}
j(\mathcal E_{a,b,c,t})=
\frac{
\left(A(a,b,c,t)
+16abct\right)^3}{(abct)^2A(a,b,c,t)}.
\end{equation}
This implies that in general $\mathcal E_{a,b,c}$ has 
$6$ singular fibres, corresponding
to the values of $t$ for which $j(\mathcal E_{a,b,c,t})=\infty$.
The singular fibres together with their fibre type, in the general and all
special cases are given in Table~\ref{tab:types_of_fibres}.
As examples, we also tabulate the data for $\mathcal E_1$,
$\mathcal E_9$, $\mathcal E_{25}$ and $\mathcal E_{1,4,4}$.

\begin{lemma}
\label{lem:fibre_product_structure}
$X_{(a_1:a_2:a_3:a_4:a_5:a_6)}$ and $\mathcal E_{a_1,a_2,a_3}\times_{\PP^1}
\mathcal E_{a_4,a_5,a_6}$ are birational.
\end{lemma}
\begin{proof}
We can rewrite (\ref{def_by_two_equations}) as
\begin{eqnarray*}
X_1 + X_2 + X_3 &=& -(X_4 + X_5 + X_6),\\
\frac{a_1}{X_1} + \frac{a_2}{X_2} + 
\frac{a_3}{X_3} &=& -\left(
\frac{a_4}{X_4} + \frac{a_5}{X_5} + \frac{a_6}{X_6}\right),
\end{eqnarray*}
from  which,
introducing a new parameter $(\lambda_0:\lambda_1)\in\PP^1$,
 we obtain
\begin{eqnarray}
\label{eq_part1}
(X_1 + X_2 + X_3)
\left(\frac{a_1}{X_1} + \frac{a_2}{X_2} + \frac{a_3}{X_3}\right)
\lambda_0 &=&\lambda_1,\\ 
\label{eq_part2}
(X_4 + X_5 + X_6)\left(
\frac{a_4}{X_4} + \frac{a_5}{X_5} + \frac{a_6}{X_6}\right)
\lambda_0 &=&\lambda_1,
\end{eqnarray}
But these 
are equations
for $\mathcal E_{a_1,a_2,a_3}$ and $\mathcal E_{a_4,a_5,a_6}$, so
we have a map
$$
\varphi_T: X_\ba\cap T \rightarrow
({\mathcal E}_{a_1,a_2,a_3}\cap (\lambda_0\lambda_1\not=0))
\times_{\PP^1}
({\mathcal E}_{a_4,a_5,a_6}\cap (\lambda_0\lambda_1\not=0)),
$$
%given by
$$
(X_1:X_2:X_3:X_4:X_5:X_6) \mapsto
(X_1:X_2:X_3)\times(X_4:X_5:X_6).
$$
If $P:=(X_1:X_2:X_3)\in 
{\mathcal E}_{a_1,a_2,a_3}\cap (\lambda_0\lambda_1\not=0)$
and
$Q:=(X_4:X_5:X_6)\in 
{\mathcal E}_{a_4,a_5,a_6}\cap (\lambda_0\lambda_1\not=0)$,
then $(X_1+X_2+X_3)+\mu(X_4+X_5+X_6)=0$
for some unique $\mu\not=0$, and then 
$R:=(X_1:X_2:X_3:\mu X_4:\mu X_5:\mu X_6)\in X_\ba$ and $\varphi(R)=(P,Q)$.
Thus $\varphi_T$ defines a birational map.
\end{proof}

For example, this result implies that
 $X_1, X_9, X_{(1:1:1:1:4:4)}$ and $X_{(1:1:1:4:4:9)}$
are birational to
$\mathcal E_1\times_{\PP^1} \mathcal E_1$,
$\mathcal E_1\times_{\PP^1} \mathcal E_9$,
$\mathcal E_1\times_{\PP^1} \mathcal E_{1,4,4}$
and
$\mathcal E_9\times_{\PP^1} \mathcal E_{1,4,4}$
 respectively.

\begin{table}
$$
\begin{array}{llllccc}
\mathcal E_{a^2,b^2,c^2}& \infty & 0   & t^2   & (t-2a)^2  &(t-2b)^2 &(t-2c)^2  \\
\cline{2-7}
{\scriptsize{\begin{array}{l}a+b+c=t\not=0,\\ \#\{2a,2b,2c,t,0\}=5
\end{array}}} & I_6    & I_2 & I_1 & I_1&I_1&I_1\\
\rl{3}\mathcal E_{a^2,a^2,b^2}& \infty & 0   & b^2   & (2a-b)^2  &(2a+b)^2 &  \\
\cline{2-6}
{\scriptsize{\begin{array}{l}\#\{b,a,2a,0\}=4
\end{array}}} & I_6    & I_2 & I_2 & I_1&I_1&\\
\rl{3}\mathcal E_{a^2,a^2,a^2}& \infty & 0   & a^2   & 9a^2  & &  \\
\cline{2-5}
{}_{a\not=0}
 & I_6    & I_2 & I_3 & I_1&&\\
\rl{3}\mathcal E_{a^2,a^2,4a^2}& \infty  & 0 &4a^2   &16a^2 && \\
\cline{2-5}
{}_{a\not=0} & I_6  &III  & I_2 & I_1&&\\
\rl{3}\mathcal E_{a^2,b^2,c^2}& \infty & 0 & 4a^2  &4b^2 &4c^2 & \\
\cline{2-6}
{\scriptsize{\begin{array}{l}a+b+c=0,\\ \#\{a,b,c,0\}=4
\end{array}}} & I_6  & III  & I_1 & I_1 & I_1&\\
\end{array}
$$

$$
\begin{array}{lllllllllllllll}
\mathcal E_1 & \infty & 0   & 1   & 9  &&&& \mathcal E_{1,4,4} & \infty & 0   & 1   & 9  &25&\\
\cline{2-5}\cline{10-14}
             & I_6    & I_2 & I_3 & I_1&&          &&    & I_6    & I_2 & I_2 & I_1&I_1& \\
\rl{3}\mathcal E_9 & \infty & 0   & 1   & 9  &25 &&& \mathcal E_{25} & \infty & 0   & 9  &25&49\\
\cline{2-6}\cline{10-14}
             & I_6    & I_2 & I_1 & I_2&I_1&        &&      & I_6    & I_2 & I_1&I_2&I_1&\
\end{array}
$$
\caption{The fibre types of the singular fibres of $\mathcal E_\ba$}
\label{tab:types_of_fibres}
\end{table}

\subsection{Existence of smooth projective Calabi-Yau models}
\label{subsec:existance_of_small_res}

Now we can apply results of Schoen \cite{schoen} to 
determine whether $X_\ba$ has
a small projective resolution or not.

\begin{lemma}
\label{lem:no_small_res_from_schoen}
If $a,b,c,d,e,f\in k\setminus\{0\}$ with
$a+b+c+d+e+f=0$ satisfy (i)
$0\notin\{\pm a\pm b\pm c,\pm d \pm e \pm f\}$, (ii)
$(a+b+c)^2\notin \{(a+b-c)^2,(a-b+c)^2,(a-b-c)^2\}$, and
(iii) $\{a^2,b^2,c^2\}\not=\{d^2,e^2,f^2\}$,
then
$\mathcal
E_{a^2,b^2,c^2}\times_{\PP^1}\mathcal E_{d^2,e^2,f^2}$ has no small projective resolution.
\end{lemma}
\begin{proof}
This follows from \cite[Lemma 3.1 (iii)]{schoen}, as in this case 
$\mathcal E_{a,b,c}$ and $\mathcal E_{d,e,f}$ 
are not isogenous, and from Table~\ref{tab:types_of_fibres}, both
have singular fibres at $(a+b+c)^2$, with that for $\mathcal E_{a,b,c}$ being $I_1$
type.
\end{proof}

In Corollary~\ref{cor:singsofX1} we saw that
$X_1$ has a smooth projective Calabi-Yau model.
However, as we shall now see, this is not true for $X_{\phi(\bb)}$
in general.  The following was pointed out to us by J.~Koll\'ar.
\begin{proposition}
\label{existence_of_smooth_models}
Suppose $Y$ is a nodal threefold with trivial canonical bundle and
that there is no small resolution of all singularities.
Then $Y$ does not posess a smooth projective Calabi-Yau model.
\end{proposition}
\begin{proof}
Assume that $Z$ 
is a smooth projective model of $Y$ which is Calabi-Yau. 
We can successively blow up Weil divisors of $Y$ until 
we obtain a projective variety $Y'$ with
all Weil divisors of $Y'$  Cartier at all nodes. In other 
words $Y'$ is factorial. By assumption $Y'$ cannot be smooth.
We have a birational map  $f:Y' \dashrightarrow Z$. Since both $Y'$ and $Z$
have trivial (and hence nef) canonical bundle the map $f$
factors as a finite sequence of flops \cite[Theorem 4.9]{Ko} .  Since 
flops in dimension $3$ do not change the type of singularity
\cite[Theorem 2.4]{Ko} this gives a contradiction.
\end{proof}

\begin{corollary}
\label{cor:examples_with_no_smoothCY}
For $\bb\in T$, the nodal variety $X_{\phi(\bb)}$ does not
posess smooth projective Calabi-Yau model unless 
$\phi(\bb)=\phi((1:-1:b:-b:c))$, 
up to permutations, with $1\pm b\pm c$ for some sign choice.
\end{corollary}
\begin{proof}
Let $\phi(\bb)=(b_1^2:b_2^2:b_3^2:b_4^2:b_5^2:b_6^2)$ with $b_1+b_2+b_3+b_4+b_5+b_6=0$.
If for some permutation of these variables the condition of
Lemma~\ref{lem:no_small_res_from_schoen} is satisfied, 
$\mathcal E_{b_1^2,b_2^2,b_3^2}
\times_{\PP^1}\mathcal E_{b_4^2,b_5^2,b_6^2}$ has no smooth projective Calabi-Yau model.
Then 
Proposition~\ref{existence_of_smooth_models} and
Lemma~\ref{lem:fibre_product_structure} imply that $X_{\phi(\bb)}$ has no 
smooth projective Calabi-Yau model.

Suppose that conditions (i), (ii) and (iii) of 
Lemma~\ref{lem:no_small_res_from_schoen} are not satisfied
for any permutations, or sign changes preserving $\sum b_i=0$, of the $b_i$.

First note that 
$(a+b+c)^2\notin \{(a+b-c)^2,(a-b+c)^2,(a-b-c)^2\}
\iff abc(a+b)(b+c)(c+a)\not=0,$
so if $b_i+b_j\not=0$ for $1\le i,j\le 6$, condition (ii) can not fail.
If we had $b_i + b_j + b_k=0$ for all triples $1\le i,j,k\le 6$, then all
$b_i$ are equal, but this contradicts $\sum b_i=0$.  
If (iii) fails for all permutations, then 
all $b_i^2$ are equal, and if $b_i+b_j\not=0$, this again gives the 
contradiction that
all $b_i$ are equal.
Hence we can assume
$b_1+b_2=0$.  Now suppose that 
$b_i+b_j\not=0$ for $3\le i,j\le 6$.  Then (ii) holds for $b_3,b_4,b_5$.
Now $b_3+b_4+b_5+b_6=-(b_1+b_2)=0$, so (i) also holds, since $a_6\not=0$.
Similarly, (i) and (ii) hold for all other triples from $b_3,b_4,b_5,b_6$.
If (iii) fails for all of these, then again, all $b_i^2$ are equal,
but this is not possible with $\sum b_i=0$.
Thus we can assume $b_3+b_4=0$, and 
since $\sum b_i=0$ we also have $b_5+b_6=0$.

Thus, if (i), (ii) or (iii) always fail, this implies that
up to permutations, $\bb=(1:-1:b:-b:c)$.  If $1 + \pm b + \pm c=0$, then
after a sign change, we have $1 + b + c =0$, and then
$\phi(\bb)=(1:1:b:b:-c)$, and taking $b_i$ in this order
satisfies (i), (ii) and (iii) of Lemma~\ref{lem:no_small_res_from_schoen}.

Conversely, if $\bb=(1:-1:b:-b:c)$,
then (i), (ii) or (iii) always fail, unless $\phi(\bb)\not=\phi(\bc)$ 
for some $\bc$ not of this form.  As in 
Proposition~\ref{prop:number_of_nodes_on_subfamilies}, 
we obtain $\bc$ from $\bb$
by changing the signs of the $b_i$ in a nontrivial subset of the $b_i$ which
sum to $0$.  For $\bc$ to have a different form is equivalent to
$1\pm b\pm c=0$ for some choice of signs.

In case $\phi(\bb)=\phi((1:-1:b:-b:c))$, 
up to permutations, with $1\pm b\pm c\not=0$ for all sign choices,
by Proposition~\ref{prop:case_of_small_res} and
Remark~\ref{rem:explicit_small_res} there is a 
small resolution of $X_{\phi(\bb)}$,
by explicit construction.

\end{proof}

In particular, this means that $X_9$, $X_{25}$, 
$X_{(1:1:1:1:4:4)}$ and $X_{(1:1:1:4:4:9)}$ have no small 
projective resolution.
Later we will see that $X_1$ and $X_9$ have the same L-series.  The Tate
conjecture then says that there should be a correspondence between them.
In this case such a correspondence can not be a birational map, since
Corollaries~\ref{cor:examples_with_no_smoothCY} and
\ref{cor:singsofX1} imply the following result.

\begin{corollary}
$X_1$ and $X_{9}$ are not birational.
\end{corollary}

\subsection{Computation of $h^{12}$}

By Proposition~\ref{h12_of_Xa}, we know all the Hodge numbers
of the varieties $\widetilde X_\ba$ and $\widehat X_\ba$,
except for $h^{11}$ and $h^{12}$, for which we only know the
difference. 
Note that $\widetilde h^{12}=\hat h^{12}$, so
it is  enough to compute $\tilde h^{12}$.
We will use the fact that in our situation $h^{12}$ is a birational invariant.
This can either be
deduced from \cite[Corollary (4.12)]{Ko} or, as Batyrev 
has informed us, by
using motivic integration.

Then from Lemma~\ref{lem:fibre_product_structure} 
it is enough to compute $h^{12}$ for
the corresponding elliptic families fibre product.
First we need to know that the elliptic 
families in question are semistable, for which we need the following lemma.

\begin{lemma}
\label{lem:satisfying_condition_for_semi_stability}
Suppose that $b_1,\dots,b_6\not=1,1,1,1,1,2$ or $1,1,1,1,2,3$ (up to scaling,
permutations and sign changes).  Then,
after possible permutation,  we can assume that $b_1\pm b_2\pm b_3\not=0$ and
$b_4\pm b_5\pm b_6\not=0$, for all sign choices.  
\end{lemma}
\begin{proof}
Suppose that this is not the case.
Then for each partition $\{i_1,i_2,i_3\}$, $\{i_4,i_5,i_6\}$ of
$1,\dots,6$ into $2$ sets of size $3$,
we have either 
$(b_{i_1}\pm b_{i_2} \pm b_{i_3})=0$
or
$(b_{i_4}\pm b_{i_5} \pm b_{i_6})=0$, (or both)
(for some choice of signs).

First we will show that some of the $b_i$ must be equal, up to sign.
There are $10$ partitions, so we have ten distinct
sets $\{i_1,i_2,i_3\}$ of size $3$, with
$(b_{i_1}\pm b_{i_2} \pm b_{i_3})=0$.  
These sets have $30=3\cdot10$ entries, and so each of $1,\dots,6$ must
occur at least $5=30/6$ times, i.e., in at least $5$ of the $10$ sets. 

Take $5$ of the sets containing $1$.  
Suppose we always have $b_i\not=\pm b_j$ for $i\not=j$.
Then for each other element occuring in these $5$ sets, it must
occur exactly twice, with opposite signs; or else we would have,
e.g., $b_1 + \varepsilon 
b_2 \pm b_4=0$ and $b_1 + 
\varepsilon b_2 \pm b_5=0$ where $\varepsilon = \pm 1$, 
so $b_4=\pm b_5$.
But then when we take the sum of all $5$ equalities, we would get
$5b_1=0$, a contradition.  Hence we have $b_i=\pm b_j$ for some 
$i,j\in\{2,3,4,5,6\}$.  This is similarly the case for every other
subset of $5$ elements.  This implies that there are two 
(not necessarily disjoint) pairs of
equal elements amongst the $b_i$.  Thus the $b_i$ have the form
$a,a,a,b,c,d$ or $a,a,b,b,c,d$, up to sign.
In the first case we must have $b\pm c\pm d=0$.  From the other possible
partitions, we get (possibly after replacing $b$ with $-b$ or
$c$ with $-c$),
$$\begin{array}{lll}
            & b = \pm 2a & \text{or } a + c + d=0,\\
\text{and } & c = \pm 2a & \text{or } a + b \pm d=0,\\
\text{and } & d = \pm 2a & \text{or } a \pm c \pm b=0.
\end{array}
$$
We can not have all of the first column true, since this would imply
$a=0$ by taking the sum (up to sign), since $b\pm c\pm d=0$.
If two of the first column hold, eg, $b=\pm 2a$, and $c=\pm 2a$,
then from the last line we get $a=\pm 4a$, a contradiction.
If from the first column we just have
$d=\pm 2a$, then from the second, we get 
$b,c\in\{a,-3a\}$; since $b\pm c\pm d=0$, we must have $b=c=a$ or
$\{b,c\}=\{a,-3a\}$.
The first possibility gives that the $b_i$ are (after scaling)
given, up to sign, by $1,1,1,1,2,3$, and the second gives $1,1,1,1,1,2$.

If the whole of the last column is true, then substituing $b=\pm c\pm d$
and eliminating $a$, we get $c+d=\pm c\pm 2d$ or $\pm c$ and 
$c+d=\pm d\pm 2c$ or $\pm d$.  Since $c,d\not=0$, this gives
$2c = \pm d$ or $\pm 3d$ and 
$2d = \pm c$ or $\pm 3c$. 
But all possibilities lead to $c=d=0$, a contradiction.

Now suppose the $b_i$ have the form $a,a,b,b,c,d$.  This gives us
$$\begin{array}{lll}
            & c = \pm 2a & \text{or } d =\pm 2b,\\
\text{and } & d = \pm 2a & \text{or } c =\pm 2b.\\
\end{array}
$$
If we have one true from each column, we get $b=a$, and this gives
the previous case.  Otherwise, we can assume we have
$c=\pm d=\pm 2a$. The $b_i$ now must have the form
$a,a,b,b,2a,2a$.  Thus $b=a$ or $b=2a$.  This gives, up to sign change,
scaling and permutation, $1,1,1,1,2,2$ or $1,1,2,2,4,4$.
But in the first case we have $(1\pm1\pm1)(1\pm2\pm2)\not=0$,
and in the second $(1\pm2\pm4)(1\pm2\pm4)\not=0$, so neither of these
is possible.

Hence the only possibility for the $b_i$ such that there is no
perumtation with $b_1\pm b_2\pm b_3\not=0$ and 
$b_4\pm b_5\pm b_6\not=0$, are those given.
\end{proof}

\begin{proposition}
\label{prop:values_of_h12}
If $\mathcal F_i$ is the smallest set in 
Table~\ref{subfamiliesTable} containing
$\ba$, then $h^{12}(\widetilde X_\ba)=\dim \mathcal F_i$.
In particular, if $\ba\in\PP^5\setminus\phi(\PP^4)$ 
with all coordinates nonzero, then $h^{12}(\overline X_\ba)=5$.  
\end{proposition}
\begin{proof}
For some $b_i\not=0$ we have
 $\ba=(b_1^2:\cdots: b_6^2)\in \PP^5$.  We set
$\mathcal E^1=\mathcal E_{b_1^2,b_2^2,b_3^2}$ and 
$\mathcal E^2=\mathcal E_{b_4^2,b_5^2,b_6^2}$.
By Lemma~\ref{lem:fibre_product_structure}, 
$X_\ba$ is birational to $\mathcal E^1\times_{\PP^1}\mathcal E^2$, and
since $h^{12}$ is a birational invariant, it is enough to compute
 $h^{12}(\mathcal E^1\times_{\PP^1}\mathcal E^2)$.

Assuming $b_i$ are not $1,1,1,1,1,2$ or $1,1,1,1,2,3$, 
by Lemma~\ref{lem:satisfying_condition_for_semi_stability},
we can assume that $b_1\pm b_2\pm b_3, b_4\pm b_5\pm b_6\not=0$
for all sign choices, so $\mathcal E^1$ and $\mathcal E^2$
are semistable (see Table~\ref{tab:types_of_fibres}).  
Let $c_i(s)$ be the number of components of the fibre
${\mathcal E^i_s}$, $d=1$ if $\mathcal E^1$ and $\mathcal E^2$ 
are isogenous, and $0$ otherwise,
$
S_i:=\{s\in\PP^1 \;|\; {\mathcal E^i_s} \text{ is singular}\},
$
and
$S':=S_1\cap S_2\setminus\{0,\infty\}$.
We have $S_1=\{(b_1\pm b_2\pm b_3)^2\}\cup\{0,\infty\}$,
$c_1(s)=\#\{(\eps_2,\eps_3)\in\{\pm1\}^2|(b_1+\eps_2b_2+\eps_3b_3)^2=s\}$,
for $s\not=0,\infty$,
and $S_2$ and $c_2$ are given similarly in terms of $b_4,b_5,b_6$.
Then by  \cite[(7.4)]{schoen}, for the smooth resolution 
$\widetilde{\mathcal E^1\times_{\PP^1} \mathcal E^2}$ of
$\mathcal E^1\times_{\PP^1} \mathcal E^2$ we have
\begin{eqnarray*}
\label{eqn:schoen_h12_formula}
h^{12}(\widetilde{\mathcal E^1\times_{\PP^1} \mathcal E^2})&=&
\#(S_1\cup S_2) 
+ d - 5 + \!\!\!\sum_{s\in \PP^1\setminus 
(S_1\cap S_2)}\!\!\!\!\!\!\!(c_1(s)c_2(s) -1)\\
&=&
5 + d -\#S' -\sum_{s\in S'}(c_1(s) + c_2(s) - 2)\\
&=&
5 + d  -\sum_{s\in S'}(c_1(s) + c_2(s) - 1)\\
&=&
5 + d  -\sum_{s\in S'}c_1(s)c_2(s) - (c_1(s) -1)(c_2(s) - 1).
\end{eqnarray*}
If $(b_1^2:\cdots:b_6^2)\in \mathcal F_0$, i.e., if
$\sum\pm b_i\not=0$ for all sign choices, then $d=0$ and
$S'=\varnothing$, so 
$h^{12}(\overline X_\ba)=h^{12}(
\widetilde{\mathcal E^1\times\mathcal E^2})=5$, 
as claimed.

Now suppose $\sum b_i=0$.
Note we can take $d=0$, unless all
the $|b_i|$ are all equal, since $d=1$ only if $|b_i|$ have the form
$a,b,c,a,b,c$, but these can be rearanged as $a,a,b,c,c,b$ with
$b\not=c,a$, unless all are equal.
Set $A:=\sum_{s\in S'} (c_1(s)-1)(c_2(s)-1)$.
This is $0$ unless 
the $|b_i|$ have the form 
$a,a,b,c,c,b$, with $b\not=a,c$
or $a,a,a,a,b,b$, with $b\not=a$, or
$a,a,a,a,a,a$, in which case 
$A=1,2$ or $4$ respectively. Now $-h^{12} + 5 + d + A$ is equal to
\begin{eqnarray*}
&&
\#\{(\eps_2,\eps_3,\eps_5,\eps_6)\in\{\pm 1\}^4 \;|\;
 (b_1+\eps_2b_2+\eps_3b_3)^2 = (b_4+\eps_5b_5+\eps_6b_6)^2\}\\
&=&
\#\{(\eps_2,\dots,\eps_6)\in\{\pm 1\}^5
\;|\;
 b_1+\eps_2b_2+\eps_3b_3 + \eps_4b_4+\eps_5b_5+\eps_6b_6=0\}\\
&=&\text{number of singularities on }\overline X_\ba.
\end{eqnarray*}
The first equality follows since if
$u=b_1+\eps_2b_2+\eps_3b_3$, $v=b_4+\eps_5b_5+\eps_6b_6$, and
$u^2=v^2$, then exactly one of $u=\pm v$ holds, since we 
assume $b_1\pm b_2\pm b_3,b_4\pm b_5\pm b_6\not=0$.
Proposition~\ref{prop:number_of_nodes_on_subfamilies}
gives the second equality.
Now we have
\begin{eqnarray*}
h^{12}
&=& 5 - \text{ number of singularities of } \overline X_{\ba} + d + A \\
&=&5-\#(\text{sings of }\overline X_\ba)
+\left\{
\begin{array}{lll}
0 +1& \text{if} & \ba=(1,1,b,c,c,b),\;b\not=1,c,\\
0 +2& \text{if} & \ba=(1,1,1,1,b,b),\;b\not=1,\\  
1 +4 & \text{if} & \ba=(1,1,1,1,1,1),\\
0 &&\text{ otherwise.}
\end{array}
\right.
\end{eqnarray*}
The four cases correspond to
(1) $\mathcal F_6$ and $\mathcal F_{10}$, 
(2) $\mathcal F_{11}$ and $\mathcal F_{14}$, (3) $\mathcal F_{15}$,
and (4) everything else, respectively.
Comparing the number of singularities with the dimension of $\mathcal F_i$
in Table~\ref{subfamiliesTable} gives the result.

If the $b_i$ are $1,1,1,1,1,2$ or $1,1,1,1,2,3$, then instead of the
above, we can use van Geemen's point counting method \cite{geemen_werner}
to show that the
values of $h^{12}$ in these cases are also given as stated.
\end{proof}

We immediately have the following, which also can be quickly
deduced from \cite[Proposition 7.1]{schoen}.
\begin{corollary}
\label{cor:rigid_cases}
The varieties 
$X_1$, $X_9$, $X_{(1:1:1:1:4:4)}$ and $X_{(1:1:1:4:4:9)}$
are rigid, i.e., have $h^{12}=0$.
\end{corollary}

\section{Elliptic surfaces in $X_\ba$}
\label{sec:elliptic_surfaces}

\begin{definition}
\label{def:Hij}
Let $\ba=(a_1:\cdots:a_6)\in\PP^5$, 
$1\le i<j\le 5$, $i\not=j$, and $k<l<m$,
with $\{i,j,k,l,m\}=\{1,2,3,4,5\}$.
Suppose that
\begin{equation}
\label{eqn:cond_on_ai}
\; \prod_{n=1}^6a_n\not=0,\;a_i=a_j,\;\text{ and }
\sqrt a_k\pm \sqrt a_l\pm \sqrt a_m\pm \sqrt a_6\not=0.
\end{equation}
Let $H_{ij}$ denote
the hyperplane given in $T$ by $X_i+X_j=0$.  We define
$$E^{ij}_\ba:=\overline{(X_\ba\cap H_{ij}\cap T)}\subset X_\ba,$$
and let $\widetilde E^{ij}_\ba$ be the strict transform of $E^{ij}$ 
in $\widetilde{\overline X}_\ba$.
\end{definition}
Substituting $X_i=-X_j$ in Equation~\ref{big_subfamily} 
for $X_\ba\cap T$, gives
$$
(X_l+X_m +X_n)\left(\frac{a_l}{X_l}+\frac{a_m}{X_m}
+\frac{a_n}{X_n}\right)=a_6,
$$
This is the equation
for $\mathcal E_{a_l,a_m,a_n,a_6}$ (see (\ref{eqn:Elliptic_abct_family})),
and so $E_\ba^{ij}$ is birational to 
$\mathcal E_{a_l,a_m,a_n,a_6}\times\PP^1$.

\begin{remark}
For $\ba$ satisfying condition (\ref{eqn:cond_on_ai}),
in terms of the fibre product structure, $X_{\ba}
\cong \mathcal E_{a_i,a_i,a_6}\times_{\PP^1}\mathcal E_{a_k,a_l,a_m}$,
$E_\ba^{ij}$ corresponds to the component 
$L_{a_6}\times \mathcal E_{a_k,a_l,a_m,a_6}$ 
of the fibre over $a_6$, where $L_{a_6}$ is a component of
the $I_2$ fibre $\mathcal E_{a_i,a_i,a_6,a_6}$ of 
the family $\mathcal E_{a_i,a_i,a_6}$.
\end{remark}

\begin{remark}
\label{rem:structure_of_E_ij}
With $\ba$ as above, $\widetilde E_\ba^{ij}$ is isomorphic to 
$\mathcal E_{a_k,a_l,a_m,a_6}\times\PP^1$ blown up in the $6$ points
\begin{eqnarray*}
(1:0:0)\times(1:0),
&(0:1:0)\times(1:0),&
(0:0:1)\times(1:0),\\
(1:-1:0)\times(0:1),
&(0:1:-1)\times(0:1),&
(1:0:-1)\times(0:1).
\end{eqnarray*}
\end{remark}

\begin{lemma}
\label{lem:properties_of_Eij}
Let $a_i,a_j,a_k,a_l,a_m$ be as in Definition~\ref{def:Hij}.  
Then $E_\ba^{ij}$ is smooth and contains no singularities of $X_\ba\cap T$.
\end{lemma}
\begin{proof}
Smoothness follows from considering all possible local equations, and
the fact that $\sqrt a_k\pm \sqrt a_l\pm \sqrt a_m\pm \sqrt a_6\not=0$
for all possible sign choices, so $\mathcal E_{a_k,a_l,a_m,a_6}$
is smooth.

By Lemma~\ref{sings_on_T},
a singularity on $X_\ba$ has the form
$\bb=(b_1:\dots:b_5)$ with $b_i^2=a_i$,
and $(\sum b_i)^2=a_6$.
But if $\bb\in E_\ba^{ij}$, then $b_i+b_j=0$, 
which implies $b_k+b_l+b_m\pm\sqrt{a_6}=0$.
However this contradicts condition (\ref{eqn:cond_on_ai}).
\end{proof}

\begin{definition}
With $\ba$, $i,j,k,l,m$ as in Definition~\ref{def:Hij},
and $\mathcal E_{a_k,a_l,a_m,a_6}$ as in 
Definition~\ref{def:elliptic_curve_family},
let $\phi_{ij}$ be the birational map defined where all coordinates are
nonzero by
\begin{eqnarray*}
\phi_{ij}:\mathcal E_{a_k,a_l,a_m,a_6}\times\PP^1&\dashrightarrow &X_\ba,\\
(x:y:z)\times(r:s) &\mapsto & (X_1:X_2:X_3:X_4:X_5),\\
\text{with } (X_i:X_j:X_l:X_m:X_n)&=&(rz:-rz:sx:sy:sz).
\end{eqnarray*}
\end{definition}

We want to consider the subspace of $H_3(\widetilde{\overline X}_\ba,\ZZ)$, 
spanned by the images
of the induced maps 
$$\phi_{ij}^\ast:
H_1(\mathcal E_{a_k,a_l,a_m,a_6},\ZZ)\times H_2(\PP^1,\ZZ)
\rightarrow H_3(\widetilde{\overline X}_\ba,\ZZ).$$

\begin{definition}
\label{def:cycles_on_Eij}
For $\bc=(c_1,c_2,c_3,c_4)$, $c_i\not=0,\sum\pm c_i\not=0$,
fix $\alpha_\bc,\beta_\bc$ to be $1$-cycles, with classes
$[\alpha_\bc],[\beta_\bc]$ 
spanning $H_1(\mathcal E_{\bc},\ZZ)$, with
$\alpha_\bc.\beta_\bc=1$, and with
\begin{equation}
\label{assumption}
(1:0:0),(0:1:0),(0:0:1),(1:-1:0),(0:1:-1),(1:0:-1)\notin \alpha_\bc,\beta_\bc.
\end{equation}
\end{definition}

\begin{definition}
For $\ba$ satisfying condition (\ref{eqn:cond_on_ai}), 
and $\bc=(a_k,a_l,a_m,a_6)$,
define 
$3$-cycles on $X_\ba$ by
$$\alpha^{ij}:=\phi_{ij}(\alpha_\bc\times\PP^1),\;\;\;
\beta^{ij}:=\phi_{ij}(\beta_\bc\times\PP^1).$$
\end{definition}

\begin{lemma}
\label{lem:self_int_matrix}
For $\ba$ satisfying condition (\ref{eqn:cond_on_ai}), we have
$$
A:=
\left(
\begin{array}{cc}
(\alpha^{ij})^2 & \alpha^{ij}.\beta^{ij}\\
\beta^{ij}.\alpha^{ij} & (\beta^{ij})^2 
\end{array}
\right)
=
\left(
\begin{array}{cc}
0 & -2 \\
2 & \phantom{-}0
\end{array}
\right).
$$
\end{lemma}
\begin{proof}
We apply a general result in intersection theory
(see \cite[19.2.2]{fulton}). Namely, for a closed embedding
$i:Y\hookrightarrow X$ of compact oriented manifolds, and
$\alpha\in H_r(Y)$, $\beta\in H_s(Y)$, we have
$$i_*(\alpha)\cap i_*(\beta)=i_*(Y|_Y\cap \alpha\cap\beta).$$
We take $Y=\widetilde E^{ij}_\ba\sim_{\mathrm{bir}}\mathcal E\times \PP^1$,
where $\mathcal E=\mathcal E_{a_k,a_l,a_m,a_n}$ for appropriate indices,
and $X=\widetilde{\overline X}_\ba$.
We have
\begin{equation}
\label{eqn_Y|_Y}
Y|_Y=(K_X+ Y)|_Y=K_Y=-2 (\mathcal E\times\{\text{point}\}) + \sum E_j,
\end{equation}
where $E_j$ are the exceptional divisors in the blow up
$\widetilde E_\ba^{ij}\rightarrow \mathcal E\times\PP^1$.
The first equality of (\ref{eqn_Y|_Y}) follows from the fact that
$K_{X}=\sum Q_i$, where
the $Q_i$ are quadrics coming from blowing up nodes in $X_\ba\cap T$;
by Lemma~\ref{lem:properties_of_Eij} their restriction to $\widetilde 
E_\ba^{ij}$
is trivial. The second equality comes from the adjunction formula.
The third equality comes from the fact
that the canonical bundle of the ruled surface
$\mathcal E\times \PP^1$ has degree $-2$ on the general fibre.
Since $\alpha^{ij},\beta^{ij}$ are chosen to avoid the $E_j$ 
(Remark~\ref{rem:structure_of_E_ij}
and Definition~\ref{def:cycles_on_Eij}), we have
$$\alpha^{ij}.\beta^{ij}=
\alpha^{ij}\cap\beta^{ij}\cap Y|_Y=-2(\{\text{point}\}\times\PP^1)
\cap (\mathcal E\times\{\text{point}\})=-2.$$
 The claim
$(\alpha^{ij})^2=(\beta^{ij})^2=0$ can be seen directly from the geometry.
\end{proof}

\begin{lemma}
\label{lem:intersection_matrix_ijkl}
For $\ba=(a_1:\cdots:a_6)$ with
$a_i=a_j$, $a_k=a_l$, and $\pm a_m\pm a_6\notin\{0,2\sqrt a_i,2\sqrt a_k\}$,
where $\{i,j,k,l,m\}=\{1,2,3,4,5\}$, we have
$$
\left(
\begin{array}{cc}
\alpha^{ij}.\alpha^{kl} & \alpha^{ij}.\beta^{kl}\\
\beta^{ij}.\alpha^{kl} & \beta^{ij}.\beta^{kl}
\end{array}
\right)
=0.$$
\end{lemma}
\begin{proof}
First note that substituting $X_i+X_j=X_k+X_l=0$ in 
(\ref{big_subfamily})
gives $a_m=a_6$, but we have assumed $\sqrt a_m\pm \sqrt a_6\not=0$,
so $H_{ij}\cap H_{lm} \cap X_\ba\cap T=\varnothing$. 
Local considerations show that having picked $\alpha_\bc$ and $\beta_\bc$ to
avoid certain points means that these cycles can also not intersect in
$X_\ba\setminus T$.
\end{proof}

\begin{lemma}
\label{lem:intersection_matrix_ijk}
For $\ba=(a_1:\cdots:a_6)$ with 
$a_i=a_j=a_k$ and $\pm\sqrt a_l\pm\sqrt a_m\pm\sqrt a_6\notin\{\sqrt a_i,3
\sqrt a_i\}$,
$i<j<k$ and $\{i,j,k,l,m\}=\{1,2,3,4,5\}$, we have
$$
\left(
\begin{array}{cc}
\alpha^{ij}.\alpha^{ik} & \alpha^{ij}.\beta^{ik}\\
\beta^{ij}.\alpha^{ik} & \beta^{ij}.\beta^{ik}
\end{array}
\right)
=\text{sg}(\phi_{ij}^{-1}\circ\phi_{ik})
B,\;
\text{ where }
 B:=
\left(\begin{array}{cc}
\phantom{-}0&1\\
-1&0
\end{array}\right),
$$
and where sg is the sign
of $\phi_{ij}^{-1}\circ\phi_{ik}$ as a permutation
of coordinates of $\PP^2$.
\end{lemma}
\begin{proof}
By local considerations, one can show that the cycles do not meet in
$X_\ba\setminus T$.  In $T$, the surfaces $E_\ba^{ij}$ and
$E_\ba^{ik}$ meet in an elliptic curve, with $\alpha^{ij}$ and 
$\beta^{ij}$ restricting to the images of $\alpha$ and $\beta$ on this
curve, so up to sign the intersection matrix is $B$ as above.

The sign is determined by whether or not the map
$\phi_{km}^{-1}\circ \phi_{ij}$
preserves the orientation of the chosen cycles.
Suppose $i,j,k=1,2,4$. We have maps
\begin{eqnarray*}
\phi_{12}:
(x:y:z)\times(r:s) &\mapsto & (rz:-rz:sx:sy:sz),\\
\phi_{14}:
(x:y:z)\times(r:s) &\mapsto & (rz:sx:sy:-rz:sz).
\end{eqnarray*}
Since $E^{12}_\ba\cap E^{14}_\ba$
is given by setting $X_2=X_4$, the map
$\phi_{12}^{-1}\circ \phi_{14}$ is given by
$$(x:y:z)\mapsto (y:x:z).$$
This is an odd permutation of the coordinates,
and so $\alpha_{14}\cap E_{12}$ and 
$\beta_{14}\cap E_{12}$ have the opposite orientation to
$\alpha_{12}\cap E_{14}$ and
$\beta_{12}\cap E_{14}$, and so the intersection matrix is $-B$.
Similar considerations hold in general.
\end{proof}

\begin{corollary}
\label{cor:dim_of_elliptic_part}

(i)
If $\ba=(1:1:1:1:1:t)$, for $t\not=0,1,9$, then
the intersection matrix of $\alpha_\ba^{ij}$,
$\beta_\ba^{ij}$
for $1\le i<j\le 5$
is given by the block martix
$$
\begin{array}{ccccccccccc}
&
E^{12}&E^{13}&E^{14}&E^{15}&E^{23}&
E^{24}&E^{25}&E^{34}&E^{35}&E^{45}
\\
E^{12}&\phantom{-}A&+B&-B&+B&+B&-B&+B& \phantom{-}0& \phantom{-}0& \phantom{-}0\\
E^{13}&+B&\phantom{-}A&+B&-B&+B& \phantom{-}0& \phantom{-}0&-B&+B& \phantom{-}0\\
E^{14}&-B&+B&\phantom{-}A&+B& \phantom{-}0&+B& \phantom{-}0&-B& \phantom{-}0&+B\\
E^{15}&+B&-B&+B&\phantom{-}A& \phantom{-}0& \phantom{-}0&+B& \phantom{-}0&-B&+B\\
E^{23}&+B&+B& \phantom{-}0& \phantom{-}0&\phantom{-}A&+B&-B&+B&-B& \phantom{-}0\\
E^{24}&-B& \phantom{-}0&+B& \phantom{-}0&+B&\phantom{-}A&+B&+B& \phantom{-}0&-B\\
E^{25}&+B& \phantom{-}0& \phantom{-}0&+B&-B&+B&\phantom{-}A& \phantom{-}0&+B&-B\\
E^{34}& \phantom{-}0&-B&-B& \phantom{-}0&+B&+B& \phantom{-}0&\phantom{-}A&+B&+B\\
E^{35}& \phantom{-}0&+B& \phantom{-}0&-B&-B& \phantom{-}0&+B&+B&\phantom{-}A&+B\\
E^{45}& \phantom{-}0& \phantom{-}0&+B&+B& \phantom{-}0&-B&-B&+B&+B&\phantom{-}A
\end{array}
$$
This matrix has rank $8$.

(ii)
If $\ba=(1:1:1:t:t:t)$, for $t\not=0,1$, then
the intersection matrix of the $6$ three cycles $\alpha_\ba^{ij}$,
$\beta_\ba^{ij}$ for $1\le i<j\le 3$
is given by
$$
\left(
\begin{array}{ccc}
A&B&B\\
B&A&B\\
B&B&A
\end{array}
\right).
$$
This matrix has rank $4$.
\end{corollary}
\begin{proof}
This follows from Lemmas~\ref{lem:self_int_matrix},
\ref{lem:intersection_matrix_ijkl}
and \ref{lem:intersection_matrix_ijk}.
The conditions on $t$ in (i) and (ii) ensure that
$1\pm 1\pm 1\pm \sqrt t$ and $1\pm \sqrt t\pm \sqrt t\pm \sqrt t\not=0$,
so that the elliptic curves $E_\ba^{ij}$ are nonsingular.
\end{proof}

We now look at how much of $H^3(X_\ba)$ comes from the elliptic surfaces
$E_\ba^{ij}$.  The maps $\phi_{ij}$, defined when $a_i=a_j$, and when
the roots of the remaining coeffcients can not sum to zero,
gives us a homomorphism
$$H^3_{\scriptsize\mbox{\'et}}(X_\ba,\QQ_\ell)
\ra
\bigoplus_{i,j\text{ with }a_i=a_j \text{ and} \sum_{n\not=i,j}\pm a_n\not=0}
H^1(E^{ij}_\ba,\QQ_\ell)\times H^2(\PP^1,\QQ_\ell).$$
Let $W_\ba$ be the image of this map.  The dimension of $W_\ba$ can be
determined by computing the dimension of the corresponding intersection matrix,
as in the examples in Corollary~\ref{cor:dim_of_elliptic_part}.
Let $V_\ba$ be the kernel of this map.
We have a sequence
\begin{equation}
\label{eqn:seq_with_Va}
0\rightarrow V_\ba
\rightarrow
H^3_{\scriptsize\mbox{\'et}}(X_\ba,\QQ_\ell)
\rightarrow
W_\ba
\rightarrow
0
\end{equation}
of $\mathrm{Gal }(\overline\QQ/\QQ)$ representations. (Strictly speaking the 
$\mathrm{Gal }(\overline\QQ/\QQ)$ representation
lives on the dual spaces but by abuse of language we shall still refer to the cohomolgy groups as 
$\mathrm{Gal }(\overline\QQ/\QQ)$-modules.)
By basic linear algebra 
\begin{equation}
\label{eqn:sum_of_traces}
\trace\Frob_p|H^3_{\scriptsize\mbox{\'et}}(X_\ba,\QQ_\ell)
=\trace\Frob_p|V_\ba + \trace\Frob_p|W_\ba.
\end{equation}
Since all elliptic curves over $\QQ$ are modular, if $\ba\in\PP^5(\QQ)$, the
Galois representation on $W_\ba$ is given in terms of a weight $2$ 
modular form,
with level given by the conductor of the curve.  This modular form can 
be determined by counting points.  
The values 
$\trace\Frob_p|H^3_{\scriptsize\mbox{\'et}}(X_\ba,\QQ_\ell)$ can also
be determined by counting points, and so by subtraction we obtain the
traces of the representation $V_\ba$.  We will be most interested in
cases where $V_\ba$ is $2$-dimensional.  These are the rigid examples
of Corollary~\ref{cor:rigid_cases}, 
and the following nonrigid cases. 

\begin{corollary}
\label{cor:nonrigid_modular_cases}
For $\ba$ given by one of the following,
$$
\begin{array}{lll}
\ba &\dim \mathcal F_i\\
(1:1:1:1:1:25)\in\mathcal F_{1} & 4  \\
(1:1:1:9:9:9)\in\mathcal F_{4}   & 2  \\
(1:1:4:4:4:16)\in\mathcal F_{8}  & 1
\end{array}
$$
The semi-simplification of the Galois representation on 
$H^3_{\scriptsize\mbox{\'et}}(\widetilde {\overline X}_\ba,\QQ_\ell)$
splits into a sum of Galois representations corresponding to
elliptic curves $W_\ba$, 
and a $2$-dimensional Galois-representation, $V_\ba$.
\end{corollary}
\begin{proof}
The dimension of the pieces coming from elliptic
surfaces $E_\ba^{ij}$ for $1\le i\le j\le M$, with $M=5,3$ and $2$
respectively, is given by Corollary~\ref{cor:dim_of_elliptic_part},
parts (i), (ii), and Lemma~\ref{lem:self_int_matrix} repsectively.
These values of $\ba$ are in the indicated families $\mathcal F_i$,
defined in Table~\ref{subfamiliesTable},
and in no smaller families.  In each case we have that 
$\dim W_\ba=h^{12}(\widetilde{\overline X}_\ba)=\dim \mathcal F_i$, 
(see Proposition~\ref{prop:values_of_h12}) and so 
Galois representation $V_\ba$, which is the kernel, must have dimension $2$.
\end{proof}
\begin{remark}
In the next section we will compute the Galois representations
of the $V_\ba$  in the above corollary, and show that they correspond to
weight $4$ modular forms.
\end{remark}

\section{Computing the L-series} 

We will use the Lefschetz
fixed point theorem to compute coefficients of the L-series of 
$H^3_{\scriptsize{\mbox{\'et}}}(\widetilde{\overline X}_\ba)$.
This says that for a variety $Z$ defined over $\QQ$,
the number of rational points of $Z$ over $\FF_p$ is given by
\begin{eqnarray}
\label{eqn5}
\#Z(\FF_p)=\sum (-1)^{j}
\left(\text{ Frob}^{\ast}_p|_{H^j_{{\scriptsize\text{\'et}}}(Z)}\right).
\end{eqnarray}

The spaces 
$H^0_{\scriptsize{\mbox{\'et}}}(\widetilde{\overline X}_\ba)$
and
$H^6_{\scriptsize{\mbox{\'et}}}(\widetilde{\overline X}_\ba)$
are $1$ dimensional, and $\Frob_p$ acts trivially and by multiplication
by $p^3$ respectively.  The following result gives some information about the 
Galois action on
$H^2_{\scriptsize{\mbox{\'et}}}(\widetilde{\overline X}_\ba)$,
and (by duality) on 
$H^4_{\scriptsize{\mbox{\'et}}}(\widetilde{\overline X}_\ba)$.

\begin{proposition}
\label{prop:evals_equal_to_p}
For a prime $p$ of good reduction for $\widetilde{\overline X}_\ba$,
all eigenvalues of the Frobenius action of $\Frob_p$ on
$H^2_{\scriptsize{\mbox{\'et}}}(\widetilde X_\ba)$ are equal to $p$, provided
the rulings of the quadrics $Q_i$ which are obtained by blowing up the $30$ 
singularities are defined over the field ${\mathbb F}_p$. 
\end{proposition}
\begin{proof}
We claim that $H^2_{\scriptsize\mbox{\'et}}(\widetilde P)\cong 
H^2_{\scriptsize\mbox{\'et}}(X_\ba)$, where $\widetilde P$ is the ambient
toric variety.
This suffices since $\widetilde P$ is a smooth toric variety and hence
$H^2(\widetilde P,\ZZ)$ is spanned by divisors defined over $\ZZ$.
In order to prove that the restriction $H^2(\widetilde P,\QQ)\rightarrow
H^2(X_\ba,\QQ)$ is an isomorphism we proceed as follows.  
We first observe that as in
the proof of \cite[(1.28)]{clemens} one has an isomorphism 
$H^2(X_\ba,\QQ)\cong H^2(X_u,\QQ)$ 
where $X_u$ is a general (smooth) Calabi-Yau in $|-K_{\widetilde P}|$.
Using this isomorphism it is enough to show that
$H^2(\widetilde P,\QQ)\rightarrow H^2(X_u,\QQ)$ is an isomorphism.
Both vector spaces have dimension 
$26$.  In the case of $X_u$ this was shown in 
Proposition~\ref{prop:values_of_hijXu}, and for
$\widetilde P$ this is standard toric geometry. 
(Viewing $\widetilde P$ as a repeated blow up of $\PP^4$ one sees that the
Picard group of $\widetilde P$ is spanned by the pullback of the hyperplane
sections and the $5+10+10=25$ exceptional divisors.
Note also that there are $5$ $T$-invariant hyperplanes in
$\widetilde P$ which are, of course, linearly equivalent).
It is shown in the proof of \cite[Theorem 4.42]{Ba} that the Picard
group of $X_u$ is spanned by components of divisors of the form
$Y=H\cap X_u$ where $H$ is a $T$-invariant divisor on $\widetilde P$.
In our case there are no proper faces of $\Delta$ which have interior points.
Again from the proof of 
\cite[Theorem 4.42, p. 520]{Ba} one can conclude that $Y$ is always 
irreducible and
hence
$H^2(\widetilde P,\QQ)\rightarrow H^2(X_u,\QQ)$ is an epimorphism.
Since both vector spaces have the same dimension this is indeed an isomorphism.

In order to go from $X_\ba$ to $\widetilde X_\ba$ we consider the 
Leray spectral sequence
$$
0\rightarrow H^2_{\scriptsize\mbox{\'et}}(X_\ba)\rightarrow
H^2_{\scriptsize\mbox{\'et}}(\widetilde X_\ba)\rightarrow
\bigoplus\limits^{s}_{i=1}H^2_{\scriptsize\mbox{\'et}}(Q_i)
$$
Hence it is enough to check that the rulings of the quadrics $Q_i$ are
defined over $\FF_p$. 
\end{proof}

We shall see later that to compute values of $\Frob_p$ acting on
$H^3_{\scriptsize{\mbox{\'et}}}(\widetilde{\overline X}_\ba)$,
it will be enough to count points on $X_\ba$ over finite fields.

\subsection{Counting points on $\widetilde X_\ba$}

In this section we give a formula for counting
points on $\widetilde {\overline X}_\ba$ over finite fields.  
First we determine the primes of bad reduction for $X_\ba$,
since we will only count points on $X_\ba$ at primes of good reduction.

\begin{lemma}\label{lem:primes_of_bad_reduction}
Let $\ba=(a_1:\dots:a_6)\in\PP^5(\ZZ)$, with $\prod_{i=1}^5a_i\not=0$,
and let $F(\ba)$ be the degree $16$ polynomial in $\ZZ[a_1,\dots,a_6]$ given by
\begin{equation}
\label{eqn:Fba}
F(\ba)=\prod_{(\epsilon_1,\dots,\epsilon_5)\in \{{\pm1}\}^5}
\left(\sum_{i=1}^5 {\epsilon_i}\sqrt a_i + \sqrt a_6\right).
\end{equation}
If $\ba\notin\phi(\PP^4)$,
then $\overline X_\ba\otimes \FF_p$ is smooth over $\overline \FF_p$
 for $p\nmid a_1a_2a_3a_4a_5a_6 F(\ba)$.

Furthermore,
\begin{center}
\begin{tabular}{lll}
$\widetilde X_1\otimes\FF_p$, $\widetilde X_9\otimes\FF_p$ 
and $\widetilde X_{(1:1:1:9:9:9)}\otimes\FF_p$
  & are smooth if & $p\not=2,3,$\\
\multicolumn{3}{l}{
$\widetilde X_{25}\otimes\FF_p$, $\widetilde X_{(1:1:1:1:4:4)}\otimes\FF_p$
and $\widetilde X_{(1:1:4:4:4:16)}\otimes\FF_p$}\\
\multicolumn{3}{r}{
are smooth if  $p\not=2,3,5,$ }\\
$\widetilde X_{(1:1:1:4:4:9)}\otimes\FF_p$ 
&is smooth if& $p\not=2,3,5,7. $\\
\end{tabular}
\end{center}
\end{lemma}
\begin{proof}
Lemmas~\ref{exactly_30_sings} and \ref{sings_on_T} 
describe the singularities of $X_\ba$ over any field $\FF_p$ with
$p\nmid a_i$.
The resolutions we have described over $\overline\QQ$
remain resolutions over $\overline\FF_p$.
Thus if $\ba\notin\phi(\PP^4)$,
$\overline X_\ba\otimes\FF_p$ is smooth over $\overline \FF_p$
unless $\ba\equiv\phi(\bb)\mod p$ for some $\bb$.  This is the case only if
$F(\ba)\equiv 0\mod p$.

If $\ba=\phi(\bb)$, the primes of bad reduction are the prime factors of 
$a_i$, and of the nonzero factors of
$F(\ba)$.  E.g., if all $a_i=1$, then
$|\sum \pm \sqrt {a_i}|=0,2,4$ or $6$, so
the primes of bad reduction are $2, 3$.
Other examples are computed similarly.
\end{proof}

We now count points on $X_\ba$ by considering points on $X_\ba\cap T$
and on $X_\ba\setminus T$, and points added in the resolution
of singulatities, separately.

\begin{lemma}
If $\ba=(a_1:a_2:a_3:a_4:a_5:a_6)\in\PP^5(\FF^\times_p)$, then
\begin{equation}
\label{eqn:num_on_T}
\renewcommand\arraystretch{1.4}
\begin{array}{l}
\#(X_\ba\cap T)(\FF_p)\\
\phantom{M}={\displaystyle{\sum_{y,z,w=1}^{p-1}}}\left(
\left(\frac{\big((1+x+y+z)(\frac{a_2}{x} + \frac{a_3}{y} + \frac{a_4}{z} + a_5)
-a_1-a_6\big)^2-4a_1a_6}{p}\right)
+1
\right)\\
\phantom{MM}-2(p^2 - 3p+3) 
+\rho(a_1,a_6)\Big(\#\mathcal E_{a_2,a_3,a_5,a_4}(\FF_p)-6\Big),
\end{array}
\renewcommand\arraystretch{1}
\end{equation}
where $\rho(a_1,a_6)=p$ if $a_1\equiv a_6\mod p$, and $0$ otherwise,
$\left(\frac{x}{p}\right)$ is the Kronecker symbol, and
$\mathcal E_{a_2,a_3,a_5,a_4}(\FF_p)$ is the elliptic
curve given by (\ref{eqn:Elliptic_abct_family}).
\end{lemma}
\begin{proof}
We must compute the number of
solutions to (\ref{big_subfamily}) with all $X_i\not=0$.
Setting $X_5=1$,
$A=X_2 + X_3 + X_4 +1$ and $B=a_2/X_2 + a_3/X_3 + a_4/X_4 + a_5$,
and multiplying through by $X_1$, (\ref{big_subfamily}) becomes
\begin{equation}
\label{eqn:f_as_quadratic}
BX_1^2 + (AB + a_1 - a_6)X_1 + a_1A = 0.
\end{equation}
For fixed $X_2,X_3,X_4$, this has
$\left(\frac{d}{p}\right) + 1$ solutions, where 
$d=(AB - a_1-a_6) - 4a_1a_6$ is
the discriminant of (\ref{eqn:f_as_quadratic}).
This gives the term which is the sum in  (\ref{eqn:num_on_T}).
However, 
\begin{itemize}
\item[(i)] this sum 
counts solutions to (\ref{eqn:f_as_quadratic})
where $X_1=0$,
\item[(ii)]
if $A=B=0$, $a_1=a_6$, there are $p-\!\!1$ solutions, 
but the sum counts $1$.
\item[(iii)]
If $B=0, a_1=a_6$, $A\not=0$, there are $0$ solutions, but the
sum counts $1$.
\end{itemize}
If $a_1=a_6$, then (i) occurs exactly when $A=0$ and
$B\not=0$.  We must add
\begin{eqnarray}
\label{eqn:set1}
&-&\#\{(X_2,X_3,X_4)\in(\FF_p^\times)^3 | A=0\}\\
\label{eqn:set2}
&+&\#\{(X_2,X_3,X_4)\in(\FF_p^\times)^3 | A=B=0\}p\\
\label{eqn:set3}
&-&\#\{(X_2,X_3,X_4)\in(\FF_p^\times)^3 | B=0\}
\end{eqnarray}
to the sum in (\ref{eqn:num_on_T}).
Since $A=0, X_2X_3X_4\not=0$ is a plane with $3$ lines removed,
the set in line (\ref{eqn:set1})
has $p^2-3(p-1)$ points. The set in line (\ref{eqn:set3}) 
similarly has $p^2-3(p-1)$ points.
The equations $A=B=0$ for the
set in line (\ref{eqn:set2}) can be rearranged to give the
the equation
for $\mathcal E_{a_2,a_3,a_5,a_4}$, with all coordinates nonzero.
Thus in the case $a_1=a_6$ we obtain (\ref{eqn:num_on_T}).

The case $a_1\not=a_6$ is similar.
\end{proof}

\begin{lemma}
\label{lem:num_not_on_T}
We have
\begin{equation}
\label{eqn:num_not_on_T}
\#\big(X_\ba\setminus T\big)(\FF_p)=50p^2 + 10p + 20.
\end{equation}
\end{lemma}
\begin{proof}
The components of the decomposition
$X_\ba = \bigsqcup_{\sigma\in\widetilde\Sigma}(X_\ba\cap T_\sigma)$
are listed
in Table~\ref{table_of_decomposition}, and
illustrated in Figure~\ref{intersect_picture}.
If dim$T_\sigma<4$, then  
$\widetilde X_\ba\cap T_\sigma$ is rational, and so
the number of points $\#(\widetilde X_\ba\cap T_\sigma)(\FF_p)$, 
can easily be computed, and is
given in Table~\ref{counting_points_on_pieces_table}.
From Tables~\ref{table_of_decomposition} and \ref{counting_points_on_pieces_table} we have 
\begin{eqnarray*}
\#\big(X_\ba\setminus T\big)(\FF_p)&=&
  10(p^2 - 3p + 3)  +  20(2p^2-6p+5)\\
 && +  40  (p-2)  +  60(p-1)  +  30  (2p-3)  +  120,
\end{eqnarray*}
which sums to give the required result.
\end{proof}

\begin{table}
\begin{tabular}{lll}
$\widetilde {\overline X}_\ba 
=$& &The open threefold $\widetilde X_\ba\cap T$,\\
&+& $10$ translates of the surface $\widetilde X_\ba \cap (x=0)$,\\
&+& $20$ translates of the surface $\widetilde X_\ba \cap (y=0)$,\\
&+& $40$ translates of the curve  $\widetilde X_\ba \cap (x=y=0)$,\\
&+& $60$ translates of the curve  $\widetilde X_\ba \cap (x=z=0)$,\\
&+& $30$ translates of the curve  $\widetilde X_\ba \cap (y=z=0)$,\\
&+& $120$ translates of the  point $(x,y,z,w)=(0,0,0,-1)$,\\
&+& $30$ $\PP^1$s obtained in resolving the singularities in 
$X_\ba\setminus T$,\\
&+& the $Q_\bb$ obtained in the resolution of singularities in 
$X_\ba\cap T$.
\end{tabular}
\caption{Decomposition of $\widetilde {\overline X}_\ba$.
Translates mean images under the $S_5$ action.}
\label{table_of_decomposition}
\end{table}

\begin{table}
$$
\begin{array}{llll}
\multicolumn{3}{c}{\text{defining equations of }
X_\ba\cap T_\sigma
}  &\#(X_\ba\cap T_\sigma)(\FF_p)\\
x=0    & yzw\not=0 & a_1(1 + w + wz + wzy) =0   & p^2 - 3p + 3\\
y=0    & xzw\not=0 & (a_1 + a_2x)(1 + w + wz)=0 & 2p^2 -6p + 5\\
x=y=0  & zw\not=0  & a_1(1 + w + wz) =0         & p-2\\
x=z=0  & yw\not=0  & a_1(1 + w) =0              & p-1\\
y=z=0  & xw\not=0  & (a_1 + a_2x)(1 + w) =0     & 2p-3
\end{array}
$$
\caption{Number of points on $X_\ba\cap T_\sigma$}
\label{counting_points_on_pieces_table}
\end{table}

\begin{lemma}
\label{lem:formula_to_count_points}
For $\ba=(a_1:\cdots:a_6)\in\PP^5(\QQ)$, if
the big resolution $\widetilde {\overline X}_\ba$
of $\overline X_\ba$ has smooth reduction mod $p$,
then
\begin{eqnarray*}
\# \widetilde {\overline X}_\ba(\FF_p)
& = &  48p^2 + 46p + 14\\
&&\hspace{-1in}+
{\displaystyle{\sum_{x,y,z=1}^{p-1}}}\left(
\left(\frac{\big((1+x+y+z)(\frac{a_2}{x} + \frac{a_3}{y} + \frac{a_4}{z} + a_5)
-a_1-a_6\big)^2-4a_1a_6}{p}\right)
+1
\right)\\
&&
+\sum_{\bb=(b_1:\cdots:b_5)\in\PP^4,\; \phi(\bb)=\ba}
p\left(p + 1 + {\left(\frac{b_1b_2b_3b_4b_5(\sum b_i)}{p}
\right)}\!\!\!\right)\\
&&+\rho(a_1,a_6)\Big(\#\mathcal E_{a_2,a_3,a_5,a_4}(\FF_p)-6\Big),\\
\end{eqnarray*}
where $\rho(a_1,a_6)=p$ if $a_1\equiv a_6\mod p$, and $0$ otherwise,
$\left(\frac{x}{p}\right)$ is the Kronecker symbol, and
$\mathcal E_{a_2,a_3,a_5,a_4}(\FF_p)$ is the elliptic
curve given by (\ref{eqn:Elliptic_abct_family}).
\end{lemma}
\begin{proof}
Table~\ref{table_of_decomposition} lists the components of
$\widetilde{\overline X}_\ba$.
By Proposition~\ref{exactly_30_sings}, 
if $\ba\in\PP^5(\ZZ)$, all $30$ nodes on
$X_\ba\setminus T$
are defined over $\QQ$, and the $\PP^1$s added in the small resolution
 contribute $30p$ to the sum.
If $\ba=\phi(\bb)$,
Corollary~\ref{cor:nodes_of_X} gives the number of points on the quadric
$Q_\bb$ introduced when $\bb$ is blown up.
Thus the number of points
added in resolving the singularities of $X_\ba$ is
\begin{eqnarray}
\label{eqn:num_from_sings}
30p
+\sum_{\bb=(b_1:\cdots:b_5)\in\PP^4,\; \phi(\bb)=\ba}
\left(p + 1 + {\left(\frac{{b_1b_2b_3b_4b_5\sum_{i=1}^5 b_i}}
{p}\right)}\!\!\!\right)p.
\end{eqnarray}
Adding up 
(\ref{eqn:num_on_T}), (\ref{eqn:num_not_on_T}) and
(\ref{eqn:num_from_sings}) gives the result.
\end{proof}

\begin{figure}
\setlength{\unitlength}{2763sp}%
\begin{picture}(7490,3520)(1089,-2959)
\thinlines
\qbezier(8191,-991)(8056,-1471)(7066,-1651)
\qbezier(7066,-1651)(6286,-1471)(6466,-991)
\qbezier(6466,-991)(7441,-796)(8191,-991)
\qbezier(3961,-1531)(3721,-2101)(3016,-2476)
\qbezier(3016,-2476)(2461,-2476)(2461,-1966)
\qbezier(2461,-1966)(3161,-1366)(3961,-1531)
\qbezier(2236,179)(2236,-436)(1696,-1186)
\qbezier(1696,-1186)(1186,-1336)(1291,-751)
\qbezier(1291,-751)(1771,-76)(2236,179)
\put(4126,-436){\line(-1,-3){300}}
\put(3826,-1336){\line( 4,-5){300}}
\put(4126,-1711){\line( 1, 4){225}}
\put(4351,-811){\line(-3, 5){225}}
\put(1351,-1561){\line( 1, 1){750}}
\put(2101,-811){\line( 1,-1){750}}
\put(2851,-1561){\line(-1,-1){750}}
\put(2101,-2311){\line(-1, 1){750}}
\put(2101,-811){\line( 1, 3){300}}
\put(3376,239){\line(-5, 2){375}}
\put(3001,389){\line(-6,-1){900}}
\put(2101,239){\line( 2,-1){300}}
\put(2416, 74){\line( 5, 1){975}}
\put(3391,269){\line( 0,-1){ 15}}
\put(4141,-436){\line(-6, 5){858}}
\put(3841,-1306){\line(-4,-1){992}}
\put(1101,-1171){\line( 2,-3){250}}
\put(1396,-466){\line( 1, 1){705}}
\put(2116,-2311){\line( 2,-1){540}}
\put(2656,-2581){\line( 3, 1){735}}
\put(3391,-2336){\line( 6, 5){750}}
\put(1396,-481){\line(-2,-5){270}}
\put(3136,314){\line(-6,-1){888}}
\put(1711,-1231){\line( 1,-1){720}}
\put(3991,-1531){\line( 1, 3){290}}
\put(6882,-803){\line( 1, 0){1129}}
\put(8011,-803){\line( 6,-5){522}}
\put(8533,-1238){\line(-2,-1){868.400}}
\put(7664,-1671){\line(-1, 0){1130}}
\put(6534,-1671){\line(-6, 5){520.426}}
\put(6013,-1238){\line( 2, 1){869.200}}
\put(6013,-1238){\line( 0,-1){694}}
\put(6534,-1671){\line( 0,-1){695}}
\put(7664,-1671){\line( 0,-1){695}}
\put(8533,-1238){\line( 0,-1){694}}
\put(8533,-1932){\line(-2,-1){868.800}}
\put(7664,-2366){\line(-1, 0){1130}}
\put(6534,-2366){\line(-6, 5){520.918}}
\put(6016,-1591){\line( 6,-5){528}}
\put(6541,-2011){\line( 1, 0){1140}}
\put(7681,-2026){\line( 2, 1){886}}
\put(7066,-2386){\line( 0, 1){735}}
\put(8371,-2251){\vector(-2, 3){180}}
\put(6701,-2521){\vector(1, 4){100}}
\put(7201,-2521){\vector(-1, 4){100}}
\put(2641,-526){$w=0$}
\put(1450,-766){${y\!=\!0}$}
\put(2000,-1426){$z=0$}
\put(1906,609){$x=0$}
\put(6676,-526){$y=0$}
\put(6751,-1276){$w=0$}
\put(6026,-2701){$z=0$}
\put(7000,-2851){\parbox{3cm}{\small{singularity of $X_\ba$
at $(-a_1/a_2,0,0,-1)$}}}
\put(7060,-2011){\circle*{100}}
\put(8101,-2431){$x=0$}
\put(5600,-2651){\vector(1,1){780}}
\put(4000,-2551){\parbox{2.7cm}{\small{The surface\\
$S^{-45}\subset\overline{(y=0)}$\\
of Definition~\ref{def:surfacesSij}}}}
\put(4576,429){\parbox{6cm}{\small
These figures indicate how $\widetilde X_\ba$ intersects the toric orbits
$T_\sigma$, for
$T_\sigma\subset \overline{\{x=0\}}$ and $\overline{\{y=0\}}$
}}
\end{picture}
\caption{How $\widetilde X_\ba$ intersects $P\setminus T$.}
\label{intersect_picture}
\end{figure}

\subsection{Applying Livn\'e's method}

We now want to prove that the $2$ dimensional 
Galois representations $V_\ba$, for
$\ba=(1:1:1:1:1:1)$,
$(1:1:1:1:1:9)$,
$(1:1:1:1:1:25)$,
$(1:1:1:1:4:4)$,
$(1:1:1:4:4:9)$,
$(1:1:4:4:4:16)$
and $(1:1:1:9:9:9)$
(see Corollaries~\ref{cor:rigid_cases}, 
and \ref{cor:nonrigid_modular_cases})
are modular, by comapring explicit computation and comparison of the
coefficients of these L-series with coefficients of certain modular forms.
In order to do this, we apply Falting's method, as given by Serre and
Livn\'e \cite[Theorem 4.3]{Livne}.  A simplified form is as follows:

\begin{theorem}[Faltings-Serre-Livn\'e]
Let $\rho_1$ and $\rho_2$ be two $2$-adic $2$-dimensional Galois 
representations, unramified outside a set of primes $S$.  Let
$K_S$ be the smallest field containing all quadratic extensions of
$\QQ$ ramified at primes in $S$, and
let $T$ be a set of primes disjoint from $S$.  Then if
\begin{itemize}
\item[(L.1)]
$\trace\rho_1 \equiv\trace\rho_2\equiv0$
and
$\det\rho_1 \equiv\det\rho_2$,
\item[(L.2)]
$\{\Frob_p|_{K_s}\; : \; p\in T\}$ is ``non-cubic'' 
in ${\mathrm Gal}(K_s/K)$; in particular, it is sufficient for these
sets to be equal,
\item[(L.3)]for all $p\in T$,
\begin{eqnarray*}
(a)\; \trace\rho_1\Frob_p &=& \trace\rho_2\Frob_p,\\
(b)\; \det\rho_1\Frob_p &=& \det\rho_2\Frob_p,
\end{eqnarray*}
\end{itemize}
then $\rho_1$ and $\rho_2$ have isomorphic semi-simplifications.
\end{theorem}

We want to apply this result to the Galois representation on 
$H^3_{\scriptsize\mbox{\'et}}(X_\ba\otimes\bar{\FF}_p,\QQ_\ell))$ and the
modular Galois representaion corresponding to a cuspidal Hecke eigenform,
constructed by Serre and Deligne.  Thus we need to do the following:
\begin{itemize}
\item[(L.1]a)  Check that 
$\text{trace}(\text{Frob}_p|H^3_{\scriptsize\mbox{\'et}}
(X_\ba\otimes
\bar{\FF}_p,\QQ_\ell))$ is always even for
primes $p$ of good reduction.
\item[(L.1]b) Check that the coefficients $a_p$
of the modular forms in question are all even for primes
$p$ not dividing the level.
\item[(L.1]c) and (L.3b) Remark that the determinants of both representations
are given by $\chi^3$, where $\chi$ is the cyclotomic character.
\item[(L.2)] Determine a suitable set of primes $T_s$.
\item[(L.3]a) Compute 
$\text{trace}(\text{Frob}_p|H^3_{\scriptsize\mbox{\'et}}
(\widetilde{\overline X}_\ba\otimes\bar{\FF}_p))$ 
for all $p\in T_s$, and verify these are equal to the corresponding
coefficients of the modular forms.
\end{itemize}
We now treat each of these points in turn.

(L.1a) We first remark that the trace of Frobenius on 
$H^2_{\scriptsize{\mbox{\'et}}}(\widetilde X_\ba)$
is an integer multiple of $p$ (and hence by duality the
trace of Frobenius on $H^4_{\scriptsize{\mbox{\'et}}}(\widetilde X_\ba)$
is an integer multiple of $p^2$). This can be deduced either from
the proof of Proposition~\ref{prop:evals_equal_to_p} or directly from
the analogue of the Riemann hypothesis.
Thus the coefficients of the L-series are given by
\begin{equation}
\label{eqn:trace_on_H_3}
\text{trace}(\text{Frob}_p|H^3_{\scriptsize\mbox{\'et}}
({\widetilde{\overline X}_\ba\otimes{\bar{\FF}_p}}))
=p^3 
+ p(p+1)h + 1 - \#{\widetilde{\overline X}_\ba}(\FF_p)
\end{equation}
for some integer $h$.

\begin{lemma}
The numbers $\#\widetilde{\overline X}_\ba(\FF_p)$ and
$\text{trace}(\text{Frob}_p|H^3_{\scriptsize\mbox{\'et}}
({\widetilde{\overline X}_\ba}\otimes{\bar{\FF}_p}))$ 
are even for primes $p$ of good reduction.
\end{lemma}
\begin{proof}
Equation~(\ref{eqn:num_not_on_T}) in
Lemma~\ref{lem:num_not_on_T} implies that
$\#(\widetilde{\overline X}_\ba\setminus T)(\FF_p)$ is even.
On $X_\ba\cap T$ there is an involution $X_i\mapsto a_i/X_i$.  
The fixed points are
exactly the singularties of $X_\ba$ in $T$.  After blowing up, these
are replaced by a quadratic $Q$, with 
$\#Q(\FF_p)= (p+1)^2$ or $p^2+1$
by Corollary~\ref{cor:nodes_of_X}.  In either case this is an
even number and so 
$\widetilde{\overline X}_\ba(\FF_p)$ is even.
Now (\ref{eqn:trace_on_H_3}) implies that 
$\text{trace}(\text{Frob}_p|H^3_{\scriptsize\mbox{\'et}}
(\widetilde{\overline X}_\ba\otimes\bar{\FF}_p))$ 
is also even.
\end{proof}

(L.1b) The cuspidal Hecke eigen forms we are interested in are 
$f_{6}$, $f_{12}$, $g_{30}$,
$f_{30}$, $f_{30}'$, $f_{60}$,
and $f_{90}$, with $q$ expansions
starting as in (\ref{eqn:f6}),
(\ref{eqn:f12}),(\ref{eqn:f30}), (\ref{eqn:g30}), 
(\ref{eqn:f30p}), (\ref{eqn:f60}) and (\ref{eqn:f90}).
Arbitrarily many coefficients of the $q$-expansion of these forms may be
computed by Stein's {\sc Magma} package  \cite{magma}, \cite{stein}.

\begin{lemma}
If $f$ is a cuspidal Hecke eigen form of level $N$ coprime to
$2,3,5,7$, with $q$-expansion
$f=\sum_{n\ge 1} a_n q^n$, and if 
$a_p$ is even for primes $p$ with  $11\le p\le 37$, then 
$a_p$ is even for all primes $p\ge 11$.
\end{lemma}
\begin{proof}
This is proved by the same method as \cite[Proposition 4.10]{Livne}.
Tables of \cite{jones}, list all 
$C_3$ and $S_3$ extensions of $\QQ$ unramified outside $2,3,5,7$,
and one can easily compute to see that the Frobenius at $p$ acting on
any of these extensions has
order $3$ for at least one prime $p$ with $11\le p\le 37$.
\end{proof}

From this and the computation of enough terms of the $q$ expansions,
we immediately obtain the following result.
\begin{corollary}
For a prime $p\not=2,3,5,7$, the coefficient of $q^p$ in the
$q$-expansion of the modular forms 
given by 
$f_{6}$, $f_{12}$, $g_{30}$,
$f_{30}$, $f_{30}'$, $f_{60}$, and $f_{90}$ are all even.
\end{corollary}

(L.1c) and (L.3b). This is a well known consequence of Poincar\'e duality.
Let $V$ be a two dimensional piece of $H^3$ invariant under the Frobenius
homomorphism.  
There is a map $\wedge^2 V\rightarrow H^6$, which 
is nonzero by Poincar\'e duality.  
Since both spaces are $1$-dimensional, this is an
isomorphism.  It follows from the fact that $\Frob_p$ acts by multiplication
by $p^3$ on $H^6$ that the determinant of the action on $H^3$ is $p^3$.

(L.2) By \cite[Proposition 4.11 b]{Livne}, when $X_\ba$ has smooth
reduction for $p\not=2,3,5$, we can take
\begin{equation}
\label{eqn:t235}
T_{\{2,3,5\}}=\{7,11,13,17,19,23,29,31,41,43,53,61,71,73\}.
\end{equation}
Note that we can replace $7,11,13$ here by $103,59,37$ respectively.
This will become important later on when we shall determine the trace
of Frobenius on $H^3_{\scriptsize{\mbox{\'et}}}(\widetilde X_\ba)$.
If $X_\ba$ has bad reduction at $2,3,5,7$, then we use the following
easy lemma.
\begin{lemma}\label{lem:sufficient_p_for2357}
The elements of
$\mathrm{Gal}(\QQ[\sqrt{-1},\sqrt{2},\sqrt{3},\sqrt{5},\sqrt{7}]/\QQ)$ are 
given by the identity, together with $\mathrm{Frob}_p$ for $p$ in the set
\begin{eqnarray*}
T_{\{2,3,5,7\}}&:=&
\left\{ 11, 13, 17, 19, 23, 29, 31, 37, 41, 43, 47, 53, 59, 61, 71, 73, 79, 
83,
\right.\\
&&
\left.  101, 103, 107, 109, 113, 127, 173, 193, 211, 241, 281, 283, 311
\right\}.
\end{eqnarray*}
\end{lemma}

Note that in $T_{\{2,3,5,7\}}$ we may replace $11$ and $13$ by $179$ and
$157$ respectively.

(L.3a) Using the formula of Lemma~\ref{lem:formula_to_count_points}, 
we obtain the data in Table
~\ref{tab:numbers_of_points}.
\begin{table}
$$
\begin{array}{cccccccc}
p&\#X_{1} &
\#X_{9} &
\#X_{4,4} &
\#X_{4,4,9}&
\#X_{25} &
\#X_{9,9,9}&
\#X_{4,4,4,16}
\\
7 & 3720&    3160&    3360&    3172&  3000&  3092&   3120 \\
11& 9240&    7920&    8424&    7956&  7464&  7680&    7848  \\
13& 13080&   11260&   12036&   11368& 10500& 10940&   11088\\
17& 23400&   20340&   21420&   20112& 18540& 19464&   19920\\
19& 29640&   25840&   27480&   25840& 24720& 25352&   25416\\
23& 45120&   39600&   41904&   39840& 37560& 38796&   39144\\
29& 76560&   67860&   71604&   67584& 65100& 66408&   66984\\
31& 89400&   79480&   83376&   79528& 74664& 76760&   77880\\
41& 172200&  154980&  161820&  155172& 148884&151632 & 153744 \\
43& 193080&  174160&  181224&  174400& 167640& 170636& 172656\\
53& 320400&  291780&  303012&  292392& 281580&287112& 289512\\
61& 454440&  416620&  430788&  416500& 403884&408836& 412608\\
71  &  663840 &  612720  & 634320 &  613032 & 592944  & 603720 & 609168 \\
73& 712920& 658900 & 680700&   658684& 636180&647660& 654048\\ 
103& 1735320& 1628200& 1671168& 1627156& 1586040& 1608716& 1616976\\
59 & 418440&   383040&  397224&  383124&  367560&  375720&  378600\\
37 & 134760 & 120700 & 126804 & 121168 & 114900 & 118028 & 119808
\end{array}
$$
\caption{Number of points on $X_\ba(\FF_p)$,
where for a sequence $\bb=b_1,\dots,b_i$ of length $i<6$,
$X_\bb$ means $X_\ba$, where $\ba=(b_1:\cdots:b_i:1:\cdots:1)$}
\label{tab:numbers_of_points}
\end{table}
From these values we shall be able to compute $\trace\Frob_p|_{V{\ba}}$, where
$V_\ba$ is the two dimensional Galois representation given by
$H^3_{\scriptsize\mbox{\'et}}
(\widetilde{\overline X}_\ba\otimes\bar{\FF}_p)$
in the rigid case, or by the subrepresentation
given in Corollary~\ref{cor:nonrigid_modular_cases} 
for the three nonrigid cases.

In the nonrigid cases, the elliptic curves
$\mathcal E_{1,1,1,25}$, $\mathcal E_{1,9,9,9}$ and 
$\mathcal E_{4,4,4,16}$
have $j$-invariants
$11^3 1259^3 2^{-1} 3^{-3} 5^{-4}$,
$11^3 13^3 23^3 2^{-1}3^{-12}5^{-1}$ and
$71^3 2^{-4}3^{-3} 5^{-1}$ respectively, but they all have conductor
$30$.  They are isogenous to each other, since there is only one
weight $2$ level $30$ Hecke eigen newform,
$g_{30}=\sum b_n q^n$, with $b_p$ for primes
$p\in T_{2,3,5}$
(see (\ref{eqn:t235}))
as in the following table.
$$
\begin{array}{ccccccccccccccc}
p   & 7, 103 &11, 59 &13, 37 &17&19&23&29&31&41&43&53&61&71&73\\
b_p & -4& 0& 2& 6& -4& 0& -6& 8& -6& -4& -6& -10&0 & 2
\end{array}
$$
Now from (\ref{eqn:sum_of_traces}) and (\ref{eqn:trace_on_H_3}),
we have
$$\trace(\Frob_p|_{V_\ba})
=p^3 + 1 + p(1+p)h -
\#\widetilde{\overline X}_\ba(\FF_p)
-h^{12}(\widetilde{\overline X}_\ba)pb_p.$$
From this we see that using the data in Tables~\ref{subfamiliesTable}
and \ref{tab:numbers_of_points} we can  
compute the values of  $\trace(\Frob_p|_{V_\ba})$
provided we know the integer $h$. At this point we make use of a trick 
which to our knowledge goes back to van Geemen and Werner. Recall that 
by the analogue of the Riemann hypothesis the absolute value of the 
trace of Frobenius 
on an invariant $2$-dimensional piece of 
$H^3_{\scriptsize{\mbox{\'et}}}(\widetilde X_\ba)$ 
is bounded by $2p^{3/2}$.
If $p \geq 17$ then the value of $p(1+p)$ exceeds $4p^{3/2}$ and hence
$h$ and thus also $\trace(\Frob_p|_{V_\ba})$ can be determined by the
above fomula.
Using this observation we obtain values
as in Table~\ref{tab:values_of_trace_on_H_3}, where we use the
same indexing  convention as in Table~\ref{tab:numbers_of_points}.

\begin{table}
$$
\hspace{-0.5cm}
\begin{smallmatrix}
p   & 17&19&23&29&31&41&43&53&61&71&73&103&59&37&&\text{level}\\
V_1 & -126& 20& 168& 30& -88& 42& -52&198&-538&792&218&128&-660&
254 &&6 \\
V_9 & -126& 20& 168& 30& -88& 42&-52&198& -538&792&218&128&-660&
254 &&6 \\
V_{4,4}& 18& -100& 72& -234&-16&90&452&414&422&-360&26 &8&-684& -226
&&12\\
V_{4,4,9}& 102&20&-72&306&-136&-150&-292&-414&-418&480&434&1172
&-744&-214
&&60\\
V_{25}& 42&-76&0&6&-232&234& -412& 222& -490& 120& 746 &-560&660&
430 &&30\\
V_{9,9,9}&-66&-100&-132&90&152&438&32&-222& 902&-432&
362&-1812&-420&114&&90\\
V_{4,4,4,16}&-114&140&72&210&272&-198&-268&-78&302&-768&-478 &640& 240&
-260&&30
\end{smallmatrix}
$$
\caption{Values of $\trace(\Frob_p|_{V_\ba})$}
\label{tab:values_of_trace_on_H_3}
\end{table}
Note that in all of these cases  
$h=h^{11}(\widetilde{\overline X}_\ba)$. We conjecture that this
is always the case. To prove this, it would be enough to prove that 
 the action of the Frobenius on $H^2$ and $H^4$ 
 is multiplication by $p$ and $p^2$ respectively.
The above numbers can be verified to be coefficients of weight $4$
cuspidal Hecke eigen modular forms of levels
$6,6,12,60,30,90$ and $30$ respectively, as indicated in the last column.
  In Lemma~\ref{lem:primes_of_bad_reduction} we saw that
  $\widetilde X_{(1:1:4:4:9)}$ also has bad reduction at $7$.
  Thus in this case we have to compute the number of point over
  $\FF_p$ for additional primes as given in
  Lemma~\ref{lem:sufficient_p_for2357}, and from this, as before, we
  can compute the values of the traces on $V_\ba$.  This data is given in
  Table~\ref{tab:data_for_X_449}.
\begin{table}
$
%\begin{array}{lllllllllllll}
\begin{smallmatrix}
p &37    &     47&     53&     59&     61&     71&     79&     83&     101&
103     \\
\#\widetilde{\overline{X}}_\ba(\FF_p)
 &121168& 216696& 292392& 383124& 416500& 613032& 807688& 921000& 1546944&
1627156      \\
\trace &-214  & -72   & -414  & -744  & -418  & 480   & 1352  & -612  &
-1542  & 1172   \\
\\
\end{smallmatrix}
%\end{array}
$

$
\begin{smallmatrix}
p & 107    &     109&     113&     127&     173&     193&      211&
241&      281\\
\#\widetilde{\overline{X}}_\ba(\FF_p)& 1800888& 1896388& 2086824& 2863252&
6680856& 9063196& 11627272& 16915444& 26156796\\
\trace & 1956   & -1858  & 174    & -2068  & 1962   & -2038  & 3260    &
-1822   & -6654\\
\\
\end{smallmatrix}
$

$
\begin{smallmatrix}
p&      283&  311  & 179    & 157\\
\#\widetilde{\overline{X}}_\ba(\FF_p)& 26685544& 34931928& 7345764& 5110360 \\
\trace & -1756   & -96    &  576    & -166
\end{smallmatrix}
$
\caption{Values of $\#X_{\ba}(\FF_p)$
and $\trace(\Frob_p|_{V_\ba})$
for $\ba=(1:1:4:4:9)$
}
\label{tab:data_for_X_449}
\end{table}

%\newpage
Thus, applying \cite[Theorem 4.3]{Livne}, we have the following.
\begin{theorem}
For $\ba$ as in Corollaries~\ref{cor:rigid_cases}
and \ref{cor:nonrigid_modular_cases}
the two dimensional Galois representations $V_\ba$ given by
(\ref{eqn:seq_with_Va})
are modular, corresponding to the weight $4$ modular
forms, with coefficients and level
 indicated in Table~\ref{tab:values_of_trace_on_H_3}.
\end{theorem}

%%%%%%%%%%%%%%%%%%%%%%%%%%%%%%%%%%%%%%%%%%%%%%%%%%%%%%%%%%%%%%%%%%%%%%%%%%%%%%%

% References
%
%%%%%%%%%%%%%%%%%%%%%%%%%%%%%%%%%%%%%%%%%%%%%%%%%%%%%%%%%%%%%%%%%%%%%%%%%%%%%%%
%References:
\bibliographystyle{alpha}

\medskip

\noindent
Klaus Hulek,\\
Institut f\"ur Mathematik (C),\\
Universit\"at Hannover\\
Welfengarten 1, 30060 Hannover, Germany\\
{\tt hulek@math.uni-hannover.de}

\medskip

\noindent
Helena A. Verrill\\
Department of Mathematics\\
Louisiana State University\\
Baton Rouge, LA 70803-4918\\
{\tt verrill@math.lsu.edu}

\end{document}